\documentclass[letterpaper, 11pt]{article}

\usepackage{amsmath, amssymb}
\usepackage[T1]{fontenc}
\usepackage[left=2.5cm, right=2.5cm, top=2.5cm, bottom = 2.5cm]{geometry}
\usepackage[hidelinks]{hyperref} 
\usepackage[utf8]{inputenc}
\usepackage{makecell}
\usepackage{multirow}
\usepackage{placeins}
\usepackage{subcaption}
\usepackage{tikz}
\usepackage[numbers,sort]{natbib} 
\usetikzlibrary{decorations.pathreplacing}

\newtheorem{remark}{Remark}

\newcommand{\boldb}{\boldsymbol{b}}
\newcommand{\boldf}{\boldsymbol{f}}
\newcommand{\boldn}{\boldsymbol{n}}

\newcommand{\boldv}{\boldsymbol{v}}
\newcommand{\boldx}{\boldsymbol{x}}

\newcommand{\boldxi}{\boldsymbol{\xi}}
\newcommand{\boldeta}{\boldsymbol{\eta}}
\newcommand{\boldzero}{\boldsymbol{0}}

\newcommand{\tdx}[1]{\frac{\mathrm{d}{#1}}{\mathrm{d}x}}
\newcommand{\tdy}[1]{\frac{\mathrm{d}{#1}}{\mathrm{d}y}}
\newcommand{\dx}{\mathrm{d}x}
\newcommand{\dy}{\mathrm{d}y}

\newcommand{\tsum}{\texttt{sum}}
\newcommand*\samethanks[1][\value{footnote}]{\footnotemark[#1]}

\title{Matrix-oriented FEM formulation for stationary and time-dependent PDEs on x-normal domains}%

\author{Massimo Frittelli\thanks{University of Salento, Department of Mathematics and Physics ``E. De Giorgi'', Via per Arnesano, 73100 Lecce, Italy. Emails: \texttt{massimo.frittelli@unisalento.it}, \texttt{ivonne.sgura@unisalento.it}}, Ivonne Sgura\samethanks}

\date{}

\begin{document}

\maketitle

\begin{abstract}
When numerical solution of elliptic and parabolic partial differential equations is required to be highly accurate in space, the discrete problem usually takes the form of large-scale and sparse linear systems. In this work, as an alternative, for spatial discretization we provide a Matrix-Oriented formulation of the classical Finite Element Method, called MO-FEM, of arbitrary order $k\in\mathbb{N}$. On structured 2D domains (e.g. squares or rectangles) the discrete problem is then reformulated as a Sylvester matrix equation, that we solve by the \emph{reduced approach} in the associated spectral space. 
On a quite general class of domains, namely \emph{normal domains}, and even on special surfaces, the MO-FEM yields a \emph{multiterm Sylvester matrix equation} where the additional terms account for the geometric contribution of the domain shape.  In particular, we obtain a sequence of these equations after time discretization of parabolic problems by the IMEX Euler method.  We apply the matrix-oriented form of the Preconditioned Conjugate Gradient (MO-PCG) method to solve each multiterm Sylvester equation for MO-FEM of degree $k=1,\dots,4$ and for the lumped $\mathbb{P}_1$ case. We choose a matrix-oriented preconditioner with a single-term form that captures the spectral properties of the whole multiterm Sylvester operator. \\
For several numerical examples, we show a gain in computational time and memory occupation wrt the classical vector approach solving large sparse linear systems by a direct method or by the vector PCG with same preconditioning. As an application, we show the advantages of the MO-FEM-PCG to approximate Turing patterns with high spatial resolution in a reaction-diffusion PDE system for battery modeling.
\end{abstract}

\section*{Keywords}
Finite elements, Sylvester matrix equations, Reaction-diffusion, Turing pattern, Preconditioned Conjugate Gradients

\section*{Mathematics Subject Classification}
65F45, 65M60, 65N30

\section{Introduction}
We are interested in the discretisation of (i) elliptic PDEs of the form
\begin{equation}
\label{general_elliptic_equation}
-\Delta u + \gamma u = f(x,y), \qquad (x,y)\in \Omega \subset\mathbb{R}^2,
\end{equation}
where $\gamma \geq 0$, (ii) parabolic PDEs of the form
\begin{equation}
\label{general_parabolic_equation}
u_t - d_u \Delta u = f(x,y,t), \qquad (x,y,t) \in \Omega \times [0,T],
\end{equation}
with $d_u > 0$ being a diffusion coefficient, and (iii) reaction-diffusion systems (RDS) of the form
\begin{equation}
\label{general_RDS}
\begin{cases}
u_t - d_u\Delta u = f(u, v);\\
v_t - d_v\Delta v = g(u, v),
\end{cases}
\qquad (x,y,t) \in \Omega \times [0,T],
\end{equation}
with $d_u, d_v > 0$ being diffusion coeffcients. Problems \eqref{general_elliptic_equation}-\eqref{general_RDS} are endowed with either homogeneous Dirichlet or Neumann boundary conditions, problems \eqref{general_parabolic_equation}-\eqref{general_RDS} are endowed with suitable initial conditions.
The relevance of the PDE problems \eqref{general_elliptic_equation}-\eqref{general_RDS} is well-known, as they find numerous applications across all fields of science. We stress that the RDS \eqref{general_RDS} is the playground of Turing's theory of morphogenesis \cite{turing}, which encompasses extremely diverse applications such as biological patterning \cite{barreira2011surface}, biomembrane modelling \cite{elliott2010modeling}, tumour growth \cite{chaplain2001spatio}, metal dealloying \cite{eilks2008numerical}, financial risk management \cite{becherer2005classical}, oscillating chemical reactions \cite{vanag2004waves} and the recent applications to metal electrodeposition \cite{bozzini2013spatio} which we will consider in the present work. We focus on the approximation of Turing patterns, because, from a computational point of view, this is a challenging task since fine meshes are required in space to capture the morphological class of the pattern itself (spots, labyrinths, etc) that must be attained as steady state of the PDE dynamics for long time of integration.

Among the existing methods for the spatial discretisation of problems \eqref{general_elliptic_equation}-\eqref{general_RDS} we mention finite differences \cite{jordan1965calculus}, finite elements \cite{hughes2012finite}, spectral methods \cite{chaplain2001spatio}, kernel methods \cite{kansa1990multiquadrics} and many more. It is well known that numerical methods typically approximate the elliptic problem \eqref{general_elliptic_equation} through an algebraic system in vector form:
\begin{equation}
\label{general_elliptic_discrete}
A\boldxi = \boldb,
\end{equation}
with the vector $\boldxi$ containing the coefficients of the expansion of the numerical solution in a given discrete function basis, the matrix $A$ approximating the operator $\mathcal{L}(u) = -\Delta u + \gamma u$ and the vector $\boldb$ approximating the right-hand side of \eqref{general_elliptic_equation}. For the time-dependent problems \eqref{general_parabolic_equation}-\eqref{general_RDS}, a common general approach is the so-called method of lines (MOL), which consists of discretising the spatial variables with a spatial method of choice, thereby producing a continuous-in-time ODE system. For problem \eqref{general_parabolic_equation}, in a general setting including several spatial methods, the resulting spatially discrete problem takes the form of the following ODE system in vector form
\begin{align}
\label{general_parabolic_discrete}
&M\dot{\boldxi} + A_u\boldxi = \boldb(t),  \qquad \boldxi(0) = \boldxi_0,\qquad t \in [0,T],
\end{align}
with the vector $\boldxi = \boldxi(t)$ containing the time-dependent coefficients of the expansion of $u$, the matrix $A_u$ approximating the operator $\mathcal{L}(u) := -d_u \Delta u$, the vector $\boldb(t)$ approximating the right-hand-side of \eqref{general_parabolic_equation}, and the matrix $M$ depends on the spatial methods (e.g. the identity matrix for finite differences, mass matrix for finite elements). Similarly, the spatially discrete formulation of the RDS \eqref{general_RDS} becomes the following (possibly nonlinear) ODE system in vector form
\begin{align}
\label{general_RDS_discrete}
&\begin{cases}
M\dot{\boldxi} + A_u\boldxi = \boldb_1(\boldxi,\boldeta);\\
M\dot{\boldeta} + A_v\boldeta = \boldb_2(\boldxi,\boldeta);\\
\boldxi(0) = \boldxi_0, \quad \boldeta(0) = \boldeta_0,
\end{cases}
\qquad t \in [0,T].
\end{align}

The main computational challenge of \eqref{general_elliptic_discrete}-\eqref{general_RDS_discrete} is dimensionality. If the discrete function space has dimension $d$, the matrices $A, A_u, A_v, M$ appearing in \eqref{general_elliptic_discrete}-\eqref{general_RDS_discrete} are of size $d \times d$. There are special cases, as detailed below, where such matrices possess a general Kronecker structure with $n\in\mathbb{N}$ terms, e.g. 
\begin{equation}
\label{kronecker_structure}
A = \sum_{i=1}^n R_i \otimes L_i,
\end{equation}
with $\otimes$ denoting the Kronecker product and $R_i, L_i$ being matrices of lower dimension, e.g. $\sqrt{d} \times \sqrt{d}$, see for instance \cite{palitta2016matrix, sangalli2016isogeometric, powell2017efficient, mantzaflaris2017low, dautilia2020matrix, hao2020matrix}. In such cases, since
\begin{equation}
A\boldxi = \texttt{vec}\left(\sum_{i=1}^n L_i^T U R_i\right),
\end{equation}
where $U$ is such that $\texttt{vec}(U) = \boldxi$, problem  \eqref{general_elliptic_discrete} can be reformulated as the following linear algebraic matrix equation
\begin{equation}
\label{general_elliptic_matrix}
\sum_{i=1}^n L_i^T U R_i = B,
\end{equation}
where $B$ is such that $\texttt{vec}(B) = \boldb$. Problem \eqref{general_elliptic_matrix} is called a \emph{multiterm Sylvester equation}, see \cite{simoncini2016computational}. The solution of general multiterm Sylvester equations is mostly uncharted territory, as discussed in \cite{simoncini2016computational, shank2016efficient} and references therein.\\
A special case of \eqref{general_elliptic_matrix} worth mentioning is the two-term case $n=2$, when \eqref{general_elliptic_matrix} specialises to a generalised Sylvester equation
\begin{equation}
L_1^TUR_1 + L_2^TUR_2 = B,
\end{equation}
and closed-form algorithms are available, such as the Bartels-Stewart algorithm \cite{bartels1972solution} or its improvement proposed by Golub and others \cite{golub1979hessenberg}. If $L_i, R_i$ further fulfil suitable assumptions, even more efficient closed-form algorithms are available, based on spectral decomposition, see for instance \cite{dautilia2020matrix}, as we will also discuss in the next sections.\\
For the time-dependent problems \eqref{general_parabolic_discrete}-\eqref{general_RDS_discrete}, if the matrices $A_u, A_v, M$ possess a Kronecker decomposition similar to \eqref{kronecker_structure}, problems \eqref{general_parabolic_discrete}-\eqref{general_RDS_discrete} can be reformulated as \emph{matrix ODE systems}.  In this work,  we apply the Implicit-Explicit (IMEX) Euler scheme directly to \eqref{general_RDS_discrete} in vector form, then we will consider its MO counterpart.  This will yield a sequence of multiterm Sylvester matrix equations as detailed in Section \ref{sec:time_dependent_pdes}.

Matrix formulations of spatial methods for PDEs were successfully carried out in some notable cases. 
A class of elliptic problems with convection, posed on rectangular or parallelepypedal domains, was discretised via central finite differences in matrix-oriented form in \cite{palitta2016matrix} and the discrete problem takes the form of a multiterm Sylvester equation. The methodology was then extended to address more general polygonal domains, see \cite{hao2020matrix}. Elliptic anisotropic PDEs with stochastic terms were approximated via Galerkin method in matrix-oriented form in \cite{powell2017efficient}, the discrete problem is a multiterm Sylvester equation. 
Isogeometric analysis was successfully applied to various elliptic problems, see for instance \cite{sangalli2016isogeometric, mantzaflaris2017low,antolin2015efficient}. On square domains, the discrete problem is a generalized Sylvester equation, see \cite{sangalli2016isogeometric}. On more general domains defined through splines or NURBS, a Sylvester form can still be achieved if using suitable low-rank approximations of kernels, see \cite{mantzaflaris2017low}.  The work in \cite{dautilia2020matrix} addresses time-dependent problems, specifically the heat equation and RDSs, on rectangular domains, where the spatial discretisation is carried out via central finite differences in matrix-oriented form and the discrete problem takes the form of a two-term Sylvester equation, which lends itself to an extremely efficient numerical treatment based on spectral decomposition. 

In the present work we consider both elliptic PDEs of the form \eqref{general_elliptic_equation} and parabolic PDE problems such as the semilinear heat equation \eqref{general_parabolic_equation} and RDSs \eqref{general_RDS}, posed on a class of two-dimensional spatial domains known as \emph{normal domains}. For such PDE problems we propose a Matrix-Oriented Finite Element Method for the spatial discretisation, that we will define as MO-FEM. The proposed framework advances the existing theory on matrix-oriented spatial discretisation of PDEs in several directions, as listed below.

\begin{itemize}
\item To the best of the author's knowledge, the present work provides the first MO formulation of the finite element method for elliptic and parabolic PDEs. The proposed theory is general and applies to a large class of basis functions,  such as Lagrangian $\mathbb{P}_k$ basis functions, $k\in\mathbb{N}$, thereby covering arbitrarily high-order convergence in space. Special focus is given to the practical special case of lumped $\mathbb{P}_1$ finite elements. 

\item On rectangular domains, where Cartesian-structured mesh are immediate to construct, the discrete Laplacian takes the form of a two-term Sylvester equation, in analogy with matrix-oriented finite differences \cite{dautilia2020matrix} or isogeometric analysis \cite{sangalli2016isogeometric}. Moreover, since central finite differences are equivalent to lumped $\mathbb{P}_1$ finite elements, the matrix-oriented FD discretization of the considered PDE problems can be considered as a special case of the proposed theory. In this case we solve the two term Sylvester equation by the \emph{reduced approach} in the associated spectral space, see Section 4.

\item Thanks to a suitable coordinate transformation, the proposed theory applies to domains more general than rectangles, namely \emph{normal domains}. To the best of the authors' knowledge, the first applications of a matrix-oriented method to non-rectangular domains are (i) the work in \cite{mantzaflaris2017low}, where isogeometric analysis with suitable low rank approximation of kernels is applied to elliptic problems and (ii) the work in \cite{hao2020matrix}, where a matrix-oriented finite difference scheme with conformal mappings is applied to elliptic problems on \emph{polygonal domains}. On normal domains, the proposed approach adopts a curved mesh that matches the (possibly curved) boundary exactly. Hence, the proposed approach combines the low dimensionality of a matrix-oriented approach with the absence of geometric error. Since normal domains can be wrapped around a cylinder, the proposed method also applies to spatial domains that are special surfaces, namely cylinders with curvilinear edges. In all these cases, he discrete problem takes the form of a multiterm Sylvester equation, where the additional terms account for domain shape.\\
In this case,  for the numerical approximation of general multiterm Sylvester equations, we propose an iterative method: a matrix-oriented preconditioned conjugate gradient method (MO-PCG) that always converges for the considered PDE problems,thanks to the involved differential operators being coercive and self-adjoint, see \cite[Section 6.7]{saad2003iterative}. 

\item We provide numerical experiments that demonstrate that (i) for elliptic problems on square domains, both the reduced approach and the MO-PCG exhibit optimal spatial convergence, (ii) for elliptic and parabolic problems on $x$-normal domains, MO-PCG exhibits optimal convergence in space (and time, if combined with IMEX Euler), and (iii) both the reduced approach and the MO-PCG provide a significant gain in terms of computational time and memory storage in comparison to the standard vector form, solved both via a direct solver and vector PCG.

\item Concerning reaction-diffusion systems,  we show that the MO-FEM allows for accurate simulation of Turing patterns -obtained as asymptotic solutions- that might be prohibitive, in terms of time and memory, through standard finite elements in vector form, because fine spatial grids are required to capture the features of pattern morphology, after long-time integration, see \cite{dautilia2020matrix}.  In particular, here we solve a RDS of interest for battery modeling on some $x$-normal domains, cap and jar shaped,  and cylindrical surfaces with curvilinear boundaries. It is worth noting that the right-hand sides of such RDSs are not low-rank, then the solution of the multiterm Sylvester equations cannot be approximated by Krylov methods such those in \cite{shank2016efficient}.

\end{itemize} 
 
The paper is organized as follows. In Section \ref{sec:preliminaries}, we introduce preliminary definitions and results, we elaborate on the classes of spatial domains to be considered, and we introduce curvilinear Cartesian-structured meshes. In Section \ref{sec:matrix_formulation}, we define a general finite element method and we derive a Kronecker decomposition of the discrete differential operators, thereby considering as a practical variation lumped $\mathbb{P}_1$ finite elements.

In Section \ref{sec:poisson_square},  we introduce the MO formulation for these FEMs. on different kind of spatial 2D domains. Specifically,we discretise the stationary PDE problem \eqref{general_elliptic_equation} on a square domain and we solve the corresponding Sylvester equation by a spectral (\emph{reduced}) approach, also in the case of Lumped FEM. In Section \ref{sec:implementation-cap-shaped-domain-dirichlet-lumped}, we present the discretisation of the elliptic problem \eqref{general_elliptic_equation} on \emph{$x$-normal domains} and we present the solution of the corresponding multiterm Sylvester equations by the proposed matrix-oriented PCG method for FEM of orders $k=1,\dots,4$ in space.  A comparison with the vector PCG is provided in the numerical examples.

In Section \ref{sec:time_dependent_pdes}, we extend the proposed MO-FEM approach to the semilinear heat equation \eqref{general_parabolic_equation} and we apply the MO-PCG method to solve the sequence of Sylvester multiterm equations by the application of the IMEX-Euler method in time.  Specifically,  in Section 6.1, we present the convergence results on a cap-shaped domain and the computational performance in terms of execution time. 
In Section \ref{sec:turing_patterns} we present the numerical simulations of electrochemical patterns arising in batteries for different choice of the parameters in the reaction-diffusion system \eqref{general_RDS} corresponding to the DIB model \cite{lacitignola2017turing, sgura2019parameter} yielding Turing patterns with spots-worms and holes. We show that these patterns can be seen as PDE solutions on cylindrical surfaces. In Section \ref{sec:conclusions}, we provide some concluding remarks and highlight future research directions.

 
 \section{Normal domains and curvilinear structured meshes}
\label{sec:preliminaries}
In this section we introduce the classes of domains to be considered in this work and we construct suitable curvilinear meshes whose nodes possess a Cartesian ordering and that match curved boundaries exactly. We also introduce preliminary definitions and results that we will adopt in the construction of the proposed matrix-oriented finite element method. 

\subsection{Square domains}
To make the reader familiar with the proposed setting, we start by discretising square domains. On the one-dimensional unitary domain: $K := [0,1]$, we consider the equally spaced mesh $K_h$ with $(N+1)$ nodes, with $N\in\mathbb{N}$. For $m=0,\dots,N$, we define the $m$-th node as $x_m = m/N$. Each element of $K_h$ is of the form $E_m=[x_m,x_{m+1}]$ for some $m=0,\dots,N-1$.\\
On the square two-dimensional domain: $\Omega := [0,1]^2$ we consider the Cartesian mesh $\Omega_h$ with $(N+1)\times(N+1)$ nodes. For $m,n=0,\dots,N$, we define the node $\boldx_{mn} := (x_m, y_n) = (m/N, n/N)$. Each element of $\Omega_h$ is of the form $Q_{mn}=[x_m,x_{m+1}]\times [y_{n},y_{n+1}]$ for some $m,n=0,\dots,N-1$, see Fig. \ref{fig:square_mesh} for an illustration.\\
The above construction has a tensor structure: $\Omega = K \times K$, $\Omega_h = K_h \times K_h$ and $Q_{mn} = E_m \times E_n$.

\subsection{$x$-normal domains}
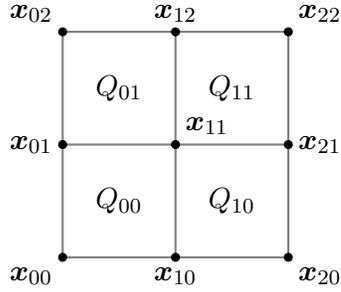
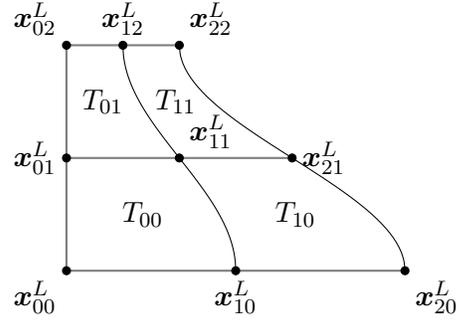
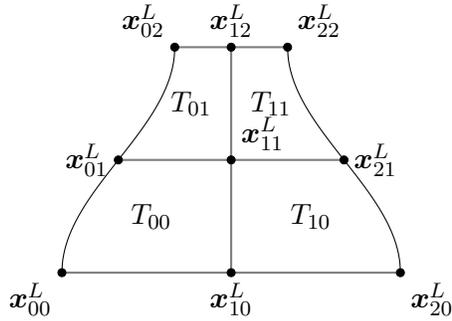
\begin{figure}
\begin{center}
\begin{subfigure}[t]{0.4\textwidth}
\begin{tikzpicture}[scale = 1.5]
\draw[gray, thick] (0,0)--(2,0)--(2,2)--(0,2)--cycle;
\draw[gray, thick] (1,0)--(1,2);
\draw[gray, thick] (0,1)--(2,1);
\draw (0.5,0.5) node[anchor = center] {$Q_{00}$};
\draw (1.5,0.5) node[anchor = center] {$Q_{10}$};
\draw (0.5,1.5) node[anchor = center] {$Q_{01}$};
\draw (1.5,1.5) node[anchor = center] {$Q_{11}$};
\filldraw[black] (0,0) circle (1pt) node[anchor=north east] {$\boldx_{00}$};
\filldraw[black] (1,0) circle (1pt) node[anchor=north] {$\boldx_{10}$};
\filldraw[black] (2,0) circle (1pt) node[anchor=north west] {$\boldx_{20}$};
\filldraw[black] (0,1) circle (1pt) node[anchor=east] {$\boldx_{01}$};
\filldraw[black] (1,1) circle (1pt) node[anchor=south west] {$\boldx_{11}$};
\filldraw[black] (2,1) circle (1pt) node[anchor=west] {$\boldx_{21}$};
\filldraw[black] (0,2) circle (1pt) node[anchor=south east] {$\boldx_{02}$};
\filldraw[black] (1,2) circle (1pt) node[anchor=south] {$\boldx_{12}$};
\filldraw[black] (2,2) circle (1pt) node[anchor=south west] {$\boldx_{22}$};
\end{tikzpicture}
\caption{Unit square $\Omega$ and its approximation $\Omega_h$ for $N=2$.}
\label{fig:square_mesh}
\end{subfigure}
\hspace*{0.05\textwidth}
\begin{subfigure}[t]{0.4\textwidth}
\begin{tikzpicture}[scale =1.5]
\draw[gray, thick] (3,0)--(0,0)--(0,2)--(1,2);
\draw[gray, thick] (0,1)--(2,1);
\draw[scale=1, domain=0:2, smooth, variable=\y, black] plot ({1+0.5*cos(90*\y)}, {\y});
\draw[scale=1, domain=0:2, smooth, variable=\y, black] plot ({2+cos(90*\y)}, {\y});
\draw ({0.5+0.25*sin(90*0.5)},0.5) node[anchor = center] {$T_{00}$};
\draw ({1.5+0.75*sin(90*0.5)},0.5) node[anchor = center] {$T_{10}$};
\draw ({0.5-0.25*sin(90*0.5)},1.5) node[anchor = center] {$T_{01}$};
\draw ({1.5-0.75*sin(90*0.5)},1.5) node[anchor = center] {$T_{11}$};
\filldraw[black] (0,0) circle (1pt) node[anchor=north east] {$\boldx_{00}^L$};
\filldraw[black] (1.5,0) circle (1pt) node[anchor=north] {$\boldx_{10}^L$};
\filldraw[black] (3,0) circle (1pt) node[anchor=north west] {$\boldx_{20}^L$};
\filldraw[black] (0,1) circle (1pt) node[anchor=east] {$\boldx_{01}^L$};
\filldraw[black] (1,1) circle (1pt) node[anchor=south west] {$\boldx_{11}^L$};
\filldraw[black] (2,1) circle (1pt) node[anchor=west] {$\boldx_{21}^L$};
\filldraw[black] (0,2) circle (1pt) node[anchor=south east] {$\boldx_{02}^L$};
\filldraw[black] (0.5,2) circle (1pt) node[anchor=south] {$\boldx_{12}^L$};
\filldraw[black] (1,2) circle (1pt) node[anchor=south west] {$\boldx_{22}^L$};
\end{tikzpicture}
\caption{An $x$-normal domain $\Omega^L$ with $L(y) = 2 + \cos 2\pi y$ and its approximation $\Omega_h^L$ for $N=2$.}
\label{fig:cap-shaped-domain_mesh}
\end{subfigure}

\vspace{10mm}

\begin{subfigure}{0.4\textwidth}
\begin{tikzpicture}[scale = 1.5]
\draw[gray, thick] (1.5,0)--(-1.5,0);
\draw[gray, thick] (-0.5,2)--(0.5,2);
\draw[gray, thick] (-1,1)--(1,1);
\draw[gray, thick] (0,0)--(0,2);
\draw[scale=1, domain=0:2, smooth, variable=\y, black] plot ({1+0.5*cos(90*\y)}, {\y});
\draw[scale=1, domain=0:2, smooth, variable=\y, black] plot ({-1-0.5*cos(90*\y)}, {\y});
\draw ({-1*sin(90*0.5)},0.5) node[anchor = center] {$T_{00}$};
\draw ({1*sin(90*0.5)},0.5) node[anchor = center] {$T_{10}$};
\draw ({-0.5*sin(90*0.5)},1.5) node[anchor = center] {$T_{01}$};
\draw ({0.5*sin(90*0.5)},1.5) node[anchor = center] {$T_{11}$};
\filldraw[black] (-1.5,0) circle (1pt) node[anchor=north east] {$\boldx_{00}^L$};
\filldraw[black] (0,0) circle (1pt) node[anchor=north] {$\boldx_{10}^L$};
\filldraw[black] (1.5,0) circle (1pt) node[anchor=north west] {$\boldx_{20}^L$};
\filldraw[black] (-1,1) circle (1pt) node[anchor=east] {$\boldx_{01}^L$};
\filldraw[black] (0,1) circle (1pt) node[anchor=south west] {$\boldx_{11}^L$};
\filldraw[black] (1,1) circle (1pt) node[anchor=west] {$\boldx_{21}^L$};
\filldraw[black] (-0.5,2) circle (1pt) node[anchor=south east] {$\boldx_{02}^L$};
\filldraw[black] (0,2) circle (1pt) node[anchor=south] {$\boldx_{12}^L$};
\filldraw[black] (0.5,2) circle (1pt) node[anchor=south west] {$\boldx_{22}^L$};
\end{tikzpicture}
\caption{A symmetric $x$-normal domain $\Omega^S$ with $S(y) = 1 + \frac{1}{2}\cos 2\pi y$ and its approximation $\Omega_h^S$ for $N=2$.}
\label{fig:symmetric_cap-shaped-domain_mesh}
\end{subfigure}
\end{center}
\caption{Pictorial illustration of the classes of spatial domains considered in this work, together with their curvilinear meshes. }
\end{figure}

We now introduce the more general class of $x$-normal domains, we discretise such domains through a curvilinear mesh and we provide some related results. Consider a smooth Cartesian curve $x = L(y)$ for $y \in [0,1]$ such that $L(y) > 0$ for all $y \in [0,1]$.
Consider the following \emph{$x$-normal domain}
\begin{equation}
\label{x_normal_domain}
\Omega^L := \{(x,y)\in\mathbb{R}^2 | y\in [0,1],\ 0 \leq x \leq L(x)\}.
\end{equation}
The reference domain $\Omega = [0,1]^2$ and the $x$-normal domain $\Omega^L$ are linked by the diffeomorphism $\eta: \Omega \rightarrow \Omega^L$ defined as follows
\begin{equation}
\label{transformation}
(x^L, y^L) = \eta(x,y) = (xL(y),y), \qquad \forall (x,y) \in \Omega.
\end{equation}
The idea of mapping a class of domains onto the reference square in order to exploit the tensor structure of the mesh is reminiscent, for example, of the work in \cite{hao2020matrix}. Let $\boldx^L_{mn} := \eta(\boldx_{mn})$, $m,n=0,\dots,N$ be the transformed nodes and let $\Omega^L_h := \eta(\Omega_h)$ be the transformed (curvilinear) mesh. For each element $Q\in\Omega_h$ let $T := \eta(Q)$ be the corresponding transformed element, see Fig. \ref{fig:cap-shaped-domain_mesh} for an illustration. The Jacobian $J_\eta$ of $\eta$ is given by
\begin{equation}
J_\eta(x,y) = \begin{pmatrix}
L(y) & xL'(y)\\
0 & 1
\end{pmatrix},
\end{equation}
so that 
\begin{equation}
\det J_\eta(x,y) = L(y).
\end{equation}
The inverse transformation $\eta^{-1}:\Omega^L \rightarrow\Omega$ is given by
\begin{equation}
(x,y) = \eta^{-1}(x^L, y^L) = \left(\frac{x^L}{L(y^L)}, y^L\right),
\end{equation}
so the Jacobian $J_{\eta^{-1}}$ of $\eta^{-1}$ is given by
\begin{equation}
J_{\eta^{-1}}(x,y) = \begin{pmatrix}
\frac{1}{L(y^L)} & -\frac{x^L L'(y)}{L^2(y)}\\
0 & 1
\end{pmatrix}.
\end{equation}
Switching back to the original coordinates, we have
\begin{equation}
J_{\eta^{-1}} \circ \eta (x,y) = \begin{pmatrix}
\frac{1}{L(y)} & -\frac{x L'(y)}{L(y)}\\
0 & 1
\end{pmatrix},
\end{equation}
which implies that
\begin{equation}
H(x,y) := (J_{\eta^{-1}} \circ \eta(x,y))(J_{\eta^{-1}} \circ \eta (x,y))^T = \begin{pmatrix}
\frac{1}{L^2(y)}+\frac{x^2 L'^2(y)}{L^2(y)} & -\frac{x L'(y)}{L(y)}\\
-\frac{x L'(y)}{L(y)} & 1
\end{pmatrix}.
\end{equation}
We finally define the matrix
\begin{equation}
\label{hat-H}
\widehat{H}(x,y) := H(x,y)\det J_\eta(x,y) = \begin{pmatrix}
\frac{1}{L(y)}+\frac{x^2 L'^2(y)}{L(y)} & -xL'(y)\\
-xL'(y) & L(y)
\end{pmatrix},
\end{equation}
which is symmetric and uniformly positive definite on $\Omega$.

\begin{remark}[Symmetric $x$-normal domains]
The proposed theory still holds true on \emph{symmetric $x$-normal domains} of the form
\begin{equation}
\label{symmetric_x_normal_domain}
\Omega^S := \{(x,y)\in\mathbb{R}^2 | y\in [0,1],\ |x| \leq S(y)\},
\end{equation}
where $S(y) > 0$, $y\in [0,1]$ is a smooth function. In this case, equations \eqref{transformation}-\eqref{hat-H} hold true by setting $\Omega = [-1/2,1/2]\times [0,1]$ and $L(y) = 2S(y)$. A symmetric $x$-normal domain $\Omega^S$ of the form \eqref{symmetric_x_normal_domain} with its discretisation $\Omega_h^S$ is shown in Fig. \ref{fig:symmetric_cap-shaped-domain_mesh}. It is worth remarking that the proposed theory can be easily extended to \emph{non-symmetric} $x$-normal domains, we do not consider this case for ease of presentation. 
\end{remark}

\begin{figure}
\begin{center}
\begin{subfigure}{0.4\textwidth}
\includegraphics[scale=0.35]{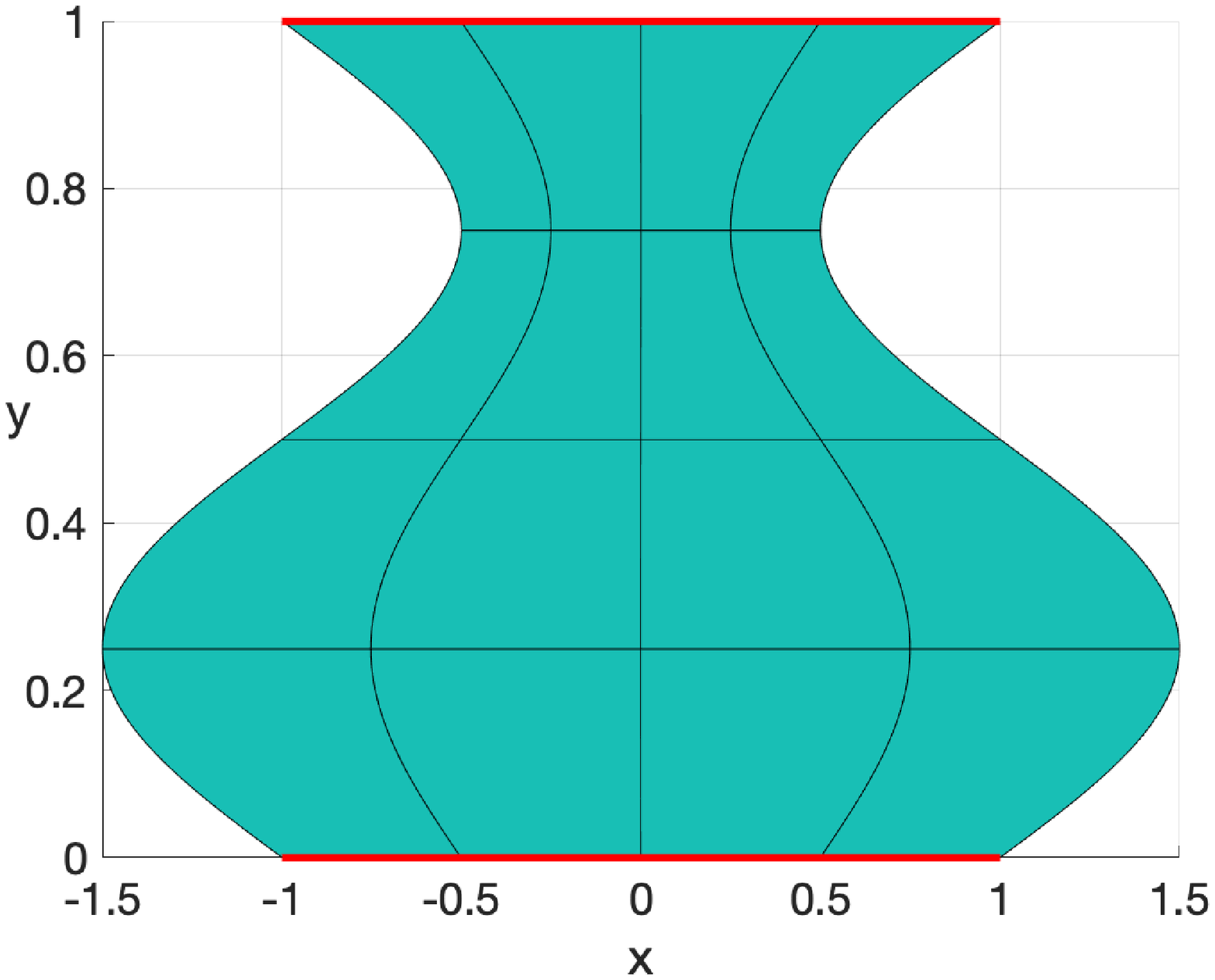}
\caption{Symmetric $x$-normal domain $\Omega^S$ of the form \eqref{symmetric_x_normal_domain} with $S(y) = 1+\frac{1}{2}\sin 2\pi y$ with its mesh $\Omega^S_h$.}
\label{fig:jar_shaped_domain}
\end{subfigure}
\hspace*{0.05\textwidth}
\begin{subfigure}{0.4\textwidth}
\includegraphics[scale=0.35]{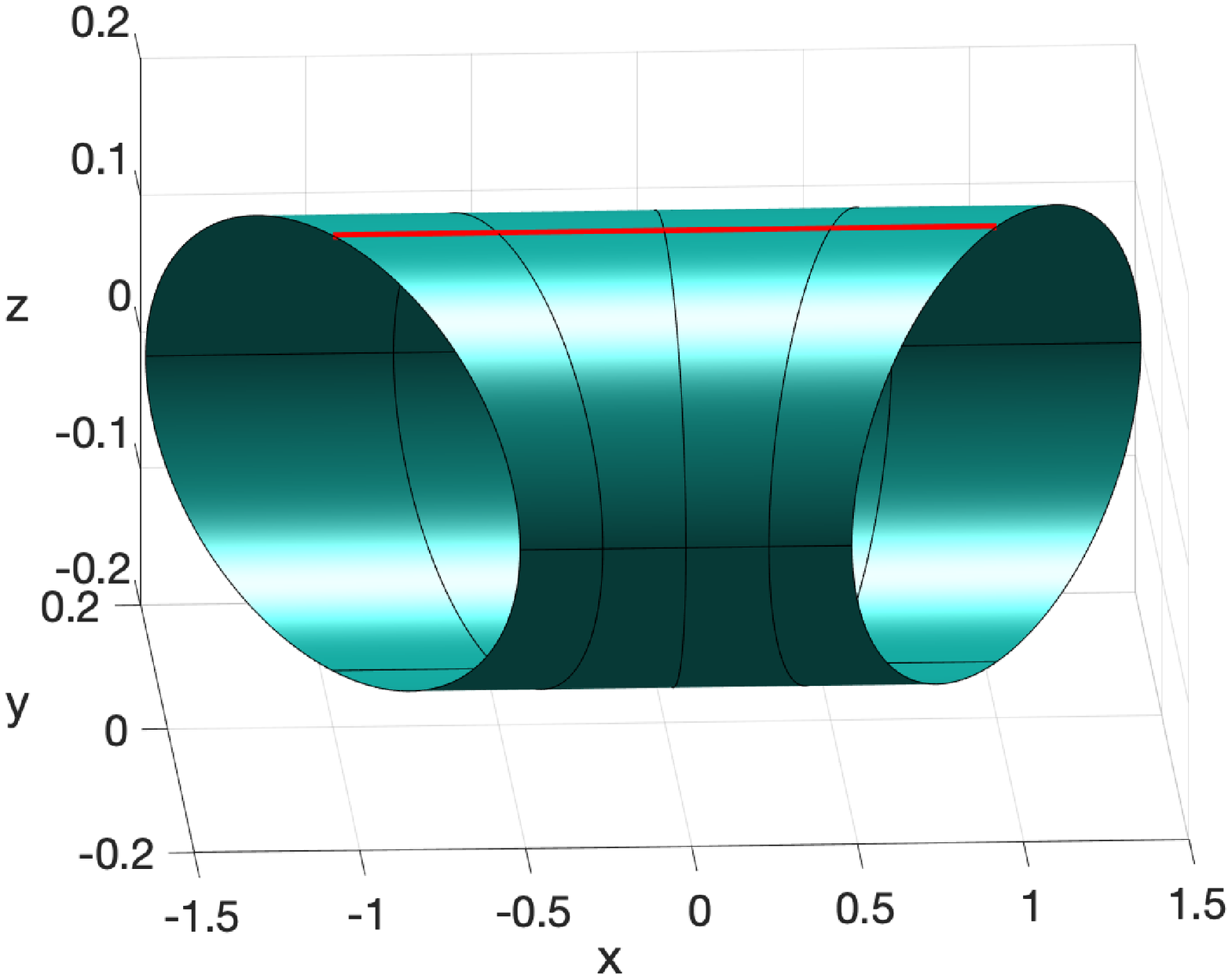}
\caption{Curvilinear cylinder $\Gamma$ defined in \eqref{definition_cylindrical_surface} with its mesh $\Gamma_h$ corresponding to the 2D x-normal domain in (a).}
\label{fig:wrapper_jar}
\end{subfigure}
\end{center}
\caption{An $x$-normal domain $\Omega^S$ and its mesh $\Omega_h^S$ are transformed into the curvilinear cylinder $\Gamma$ and its mesh $\Gamma_h$, respectively. The transformation $\sigma$ in \eqref{transformation_cylinder} joins the top- and bottom edges of $\Omega_S$, highlighted in red. }
\label{fig:cap-shaped-domain_domain_illustration}
\end{figure}

\begin{remark}[Curvilinear cylindrical surfaces]
Every $x$-normal domain \eqref{x_normal_domain} satisfying $L(0) = L(1)$ or symmetric $x$-normal domain \eqref{symmetric_x_normal_domain} satisfying $S(0) = S(1)$ can be wrapped around a cylinder through the transformation $\sigma:\Omega^L\rightarrow\mathbb{R}^3$ defined by
\begin{equation}
\label{transformation_cylinder}
\sigma(x,y) = \left(x,\frac{\sin 2\pi y}{2\pi}, \frac{\cos 2\pi y}{2\pi}\right), \qquad (x,y)\in \Omega^L.
\end{equation}
We can thus define the curvilinear cylinder $\Gamma$ and its curved mesh $\Gamma_h$ as
\begin{equation}
\label{definition_cylindrical_surface}
\Gamma := \sigma(\Omega^L); \qquad \Gamma_h := \sigma(\Omega^L_h),
\end{equation}
respectively, see Fig. \ref{fig:cap-shaped-domain_domain_illustration} for an illustration. We will show that the matrix approach proposed in the next section can be applied also to special surface PDEs, see for instance \cite{frittelli2019preserving}.
\end{remark}

\begin{remark}[More general domains]
The choice of normal domains is justified by the property that each entry of the matrix $\widehat{H}(x,y)$ defined in \eqref{hat-H} is a finite sum of separable terms. The proposed approach still applies if the coordinate transformation $\eta$ defined in \eqref{transformation} is such that the corresponding $\widehat{H}(x,y)$ retains this property.
\end{remark}

\section{Matrix-oriented formulation of FEM}
\label{sec:matrix_formulation}
We are now ready to (i) formulate a finite element method on the (possibly curvilinear) meshes introduced in the previous section and (ii) derive a Kronecker decomposition of the discrete differential operators involved in problems \eqref{general_elliptic_equation}-\eqref{general_RDS}. To make the reader familiar with the notations and results, we start with the case of square domains.

\subsection{Square domains}
Let $\{\psi_i\}_{i=1}^{N+1}$ be any finite element basis functions on the one-dimensional mesh $I_h$, e.g. piecewise Lagrange polynomials of any fixed degree. The corresponding stiffness- and mass matrices $A,M\in\mathbb{R}^{(N+1)\times (N+1)}$ in 1D are defined by
\begin{equation}
\label{stiffness_mass_1D}
a_{ij} := \int_0^1 \tdx{\psi_i}(x)\tdx{\psi_j}(x)\dx; \qquad m_{ij} := \int_0^1 \psi_i(x)\psi_j(x)\dx,
\end{equation}
for $i,j=1,\dots,N+1$, respectively. On $\Omega_h$ we choose the tensor-product local Lagrange basis $\{\phi_j\}_{j=1}^{(N+1)^2}$ defined as
\begin{equation}
\phi_{i+(N+1)j}(x,y) := \psi_i(x)\psi_j(y), \qquad \forall\ i,j=1,\dots,N+1.
\end{equation}
In general we will write $\phi_i(x,y) = \phi_{i_x}(x)\phi_{i_y}(y)$ for all $i=1,\dots,(N+1)^2$, meaning that every 2D basis function can be uniquely decomposed as a product of 1D basis functions. The stiffness- and mass matrices $\widetilde{A},\widetilde{M}\in\mathbb{R}^{(N+1)^2\times (N+1)^2}$ in 2D are defined by
\begin{equation}
\widetilde{a}_{ij} := \iint_\Omega \nabla \phi_i(x,y)\cdot \nabla \phi_j(x,y)\dx\dy; \qquad \widetilde{m}_{ij} := \iint_\Omega \phi_i(x,y)\phi_j(x,y)\dx\dy,
\end{equation}
for $i,j=1,\dots,(N+1)^2$, respectively. The stiffness- and mass matrices $\widetilde{A}$ and $\widetilde{M}$ in 2D fulfil the following relations:
\begin{equation}
\label{square_mass_decomposition}
\begin{split}
\widetilde{m}_{ij} = &\int_0^1\int_0^1 \psi_{i_x}(x)\psi_{i_y}(y)\psi_{j_x}(x)\psi_{j_y}(y)\dx\dy\\
= &\int_0^1\psi_{i_x}(x)\psi_{j_x}(x)\dx\int_0^1 \psi_{i_y}(y)\psi_{j_y}(y)\dy = m_{i_x j_x} m_{i_y j_y};
\end{split}
\end{equation}
\begin{equation}
\label{square_stiffness_decomposition}
\begin{split}
\widetilde{a}_{ij} &= \int_0^1\int_0^1 \nabla\left(\psi_{i_x}(x)\psi_{i_y}(y)\right)\cdot\nabla\left(\psi_{j_x}(x)\psi_{j_y}(y)\right)\dx\dy \\
&= \int_0^1\int_0^1 \frac{\partial}{\partial x}\left(\psi_{i_x}(x)\psi_{i_y}(y)\right)\frac{\partial}{\partial x}\left(\psi_{j_x}(x)\psi_{j_y}(y)\right)\dx\dy\\
&+ \int_0^1\int_0^1 \frac{\partial}{\partial y}\left(\psi_{i_x}(x)\psi_{i_y}(y)\right)\frac{\partial}{\partial y}\left(\psi_{j_x}(x)\psi_{j_y}(y)\right)\dx\dy\\
&= \int_0^1 \tdx{\psi_{i_x}}(x)\tdx{\psi_{j_x}}(x)\dx \int_0^1 \psi_{i_y}(y)\psi_{j_y}(y)\dy\\
&+ \int_0^1 \psi_{i_x}(x)\psi_{j_x}(x)\dx \int_0^1 \tdy{\psi_{i_y}}(y)\tdy{\psi_{j_y}}(y)\dy = a_{i_x, j_x} m_{i_y j_y} + m_{i_x, j_x} a_{i_y j_y}.
\end{split}
\end{equation}
Equations \eqref{square_mass_decomposition}-\eqref{square_stiffness_decomposition} translate to the following matrix identities:
\begin{align}
\label{kronecker_mass}
&\widetilde{M} = M\otimes M; \qquad \widetilde{A} = A\otimes M + M\otimes A,
\end{align}
where $\otimes$ denotes the Kronecker product.

\subsection{$x$-normal domains}
Let us now consider the more general case of $x$-normal and symmetric $x$-normal domains. We introduce the novel 1D matrices $B_1,B_2,C_1,C_2,M_1,M_2,M_3 \in\mathbb{R}^{(N+1)\times (N+1)}$, defined as follows:
\begin{align}
\label{newmatrix_1}
&b_{ij}^{1} := \int_0^1 \tdx{\psi_{i}}(x)\tdx{\psi_{j}}(x)L(x)\dx; 
\qquad b_{ij}^{2} := \int_0^1 \tdx{\psi_{i}}(x)\tdx{\psi_{j}}(x)x^2\dx;\\
&c_{ij}^{1} := \int_0^1 \tdx{\psi_{i}}(x)\psi_{j}(x)x\dx; 
\qquad \hspace{8mm} c_{ij}^{2} := \int_0^1 \tdx{\psi_{i}}(x)\psi_{j}(x)L'(x)\dx;\\
\label{newmatrix_7}
\begin{split}
&m_{ij}^{1} := \int_0^1 \psi_{i}(x)\psi_{j}(x)\frac{1}{L(x)}\dx; 
\qquad \hspace{3mm} m_{ij}^{2} := \int_0^1 \psi_{i}(x)\psi_{j}(x)\frac{L'^2(x)}{L(x)}\dx;\\  &m_{ij}^{3} := \int_0^1 \psi_{i}(x)\psi_{j}(x)L(x)\dx,
\end{split}
\end{align}
for $i,j=1,\dots,(N+1)$. On $\Omega_h^L$ we choose the transformed Lagrange basis $\{\phi_i^L\}_{i=1}^{(N+1)^2}$ where $\phi_i^L(x,y) := \phi_i(\eta^{-1}(x,y))$. The stiffness- and mass matrices $\widetilde{A},\widetilde{M}\in\mathbb{R}^{(N+1)^2\times (N+1)^2}$ in 2D are defined by
\begin{equation}
\label{x-normal-stiffness-mass}
\widetilde{a}_{ij} := \iint_{\Omega^L} \nabla \phi_i^L(x,y)\cdot \nabla \phi_j^L(x,y)\dx\dy; \qquad \widetilde{m}_{ij} := \iint_{\Omega^L} \phi_i^L(x,y)\phi_j^L(x,y)\dx\dy,
\end{equation}
for $i,j=1,\dots,(N+1)^2$, respectively. By using \eqref{hat-H}, the stiffness and mass matrices $\widetilde{A}$ and $\widetilde{M}$ in 2D fulfil the following relations:
\begin{equation}
\label{x-normal-decomposition-mass}
\begin{split}
\widetilde{m}_{ij} &= \iint_{\Omega^L} \phi_i^L(x,y) \phi_j^L(x,y)\dx\dy = \iint_\Omega \phi_i(x,y)\phi_j(x,y) L(y)\dx\dy\\
&= \int_0^1 \psi_{i_x}(x)\psi_{j_x}(x)\dx \int_0^1 \psi_{i_y}(y)\psi_{j_y}(y)L(y)\dy = m_{i_xj_x} m_{i_yj_y}^3;
\end{split}
\end{equation}
\begin{equation}
\label{x-normal-decomposition-stiffness}
\begin{split}
\widetilde{a}_{ij} &= \iint_{\Omega^L} \nabla\phi_i^L(x,y) \cdot \nabla \phi_j^L(x,y)\dx\dy\\
&= \iint_{\Omega^L} \left(\nabla\phi_i(\eta^{-1}(x,y))J_{\eta^{-1}}(x,y)\right) \cdot \left(\nabla \phi_j(\eta^{-1}(x,y))J_{\eta^{-1}}(x,y)\right)\dx\dy\\
&= \iint_\Omega \nabla\phi_i(x,y)^T H(x,y) \nabla\phi_j(x,y)\det J_\eta(x,y)\dx\dy\\
&= \int_0^1 \tdx{\psi_{ix}}(x)\tdx{\psi_{jx}}(x)\dx\int_0^1\psi_{i_y}(y)\psi_{j_y}(y)\frac{1}{L(y)}\dy\\
&+ \int_0^1 \psi_{ix}(x)\psi_{jx}(x)\dx\int_0^1\tdy{\psi_{i_y}}(y)\tdy{\psi_{j_y}}(y)L(y)\dy\\
&+ \int_0^1 \tdx{\psi_{ix}}(x)\tdx{\psi_{jx}}(x)x^2\dx\int_0^1\psi_{i_y}(y)\psi_{j_y}(y)\frac{L'^2(y)}{L(y)}\dy\\
&- \int_0^1 \tdx{\psi_{ix}}(x)\psi_{jx}(x)x\dx\int_0^1\psi_{i_y}(y)\tdy{\psi_{j_y}}(y)L'(y)\dy\\ &- \int_0^1 \psi_{ix}(x)\tdx{\psi_{jx}}(x)x\dx\int_0^1\tdy{\psi_{i_y}}(y)\psi_{j_y}(y)L'(y)\dy\\
&= a_{i_xj_x} m_{i_yj_y}^{1} + m_{i_xj_x} b_{i_yj_y}^{1} + b_{i_xj_x}^{2} m_{i_yj_y}^{2}  - c_{i_xj_x}^{1} c_{j_yi_y}^{2} - c_{j_xi_x}^{1} c_{i_yj_y}^{2}.
\end{split}
\end{equation}
Equations \eqref{x-normal-decomposition-mass}-\eqref{x-normal-decomposition-stiffness} translate to the following matrix identities:
\begin{align}
\label{kronecker_mass_distorted}
\widetilde{M} = M\otimes M_3; \qquad \widetilde{A} = A\otimes M_1 + M\otimes B_1 + B_2 \otimes M_2 - C_1\otimes C_2^T - C_1^T \otimes C_2,
\end{align}
where $A$ and $M$ are the standard stiffness- and mass matrices in 1D defined in \eqref{stiffness_mass_1D} and $B_1, B_2, C_1, C_2$, $ M_1, M_2, M_3, $ are novel matrices defined in \eqref{newmatrix_1}-\eqref{newmatrix_7}. Observe that on square domains, i.e. when $L(y) = 1$, we have $M_1 = M_3 = M$, $B_1 = A$, $M_2 = C_2 = 0$, so the matrix identities \eqref{kronecker_mass_distorted} encompass the special case \eqref{kronecker_mass} of square domains.

\subsection{$x$-normal domains: $\mathbb{P}_1$ elements with lumping}
The Kronecker decompositions \eqref{kronecker_mass_distorted} are absolutely general and encompass arbitrary choices of the FEM basis function, thereby including the case of Lagrangian spatial methods of any order. However, for the sake of practicality, we consider the special case of Lagrangian $\mathbb{P}_1$ finite elements \emph{with mass and stiffness lumping}, see \cite{nie1985lumped}. We show that such special case significantly simplifies the matrix identities \eqref{kronecker_mass_distorted}. This has two advantages: (i) existing finite element codes can be adapted to the proposed approach with minor modifications and (ii) the resulting numerical schemes, which take the form of matrix equations, become much easier to solve.\\
In the remainder of this section, we specialise $\{\psi_i\}_{i=1}^{N+1}$ to be the standard $\mathbb{P}_1$ Lagrangian (also called \emph{pyramidal}) basis functions. We consider the (diagonal) lumped mass matrix $M_0\in\mathbb{R}^{(N+1)\times (N+1)}$ in 1D defined by
\begin{equation}
m_{ij}^0 := \int_0^1 I_h(\psi_i(x)\psi_j(x))\dx = \delta_{ij}\int_0^1 \psi_i(x)\dx, \qquad i,j = 1,\dots, N+1,
\end{equation}
with $I_h$ being the element-wise interpolant operator \cite{nie1985lumped} and $\delta_{ij}$ being the Kronecker symbol. We also consider the (tri-diagonal) convection matrix $C \in\mathbb{R}^{(N+1)\times (N+1)}$ defined by
\begin{equation}
c_{ij} := \int_0^1 \tdx{\psi_i(x)}\psi_j(x)\dx, \qquad i,j = 1,\dots, N+1.
\end{equation}
Then we consider the (tri-diagonal) \emph{modified stiffness matrices} $A_1, A_2 \in\mathbb{R}^{(N+1)\times (N+1)}$ defined by
\begin{align}
&a_{ij}^{1} := \int_0^1 I_h\left(\tdx{\psi_{i}}(x)\tdx{\psi_{j}}(x)L(x)\right)\dx, \qquad a_{ij}^{2} := \int_0^1 I_h\left(\tdx{\psi_{i}}(x)\tdx{\psi_{j}}(x)x^2\right)\dx, 
\end{align}
for $i,j=1,\dots,N+1$. We finally define the following auxiliary diagonal matrices $D_1, D_2, D_3 \in\mathbb{R}^{(N+1)\times (N+1)}$ defined as follows
\begin{equation}
d_{ij}^1 = \delta_{ij}L(x_i), \qquad d_{ij}^2 = \delta_{ij}L'(x_i), \qquad d_{ij}^3 = \delta_{ij}x_i, \qquad i,j = 1,\dots, N+1.
\end{equation}
By combining \eqref{hat-H} and \eqref{x-normal-decomposition-stiffness} we get
\begin{equation}
\label{x-normal-stiffness-smart}
\widetilde{a}_{ij} = \iint_\Omega \left(\widehat{H}(x,y) \nabla\phi_i(x,y)\right) \cdot \nabla\phi_j(x,y)\dx\dy.
\end{equation}
The matrices defined by \eqref{x-normal-decomposition-mass} and \eqref{x-normal-stiffness-smart} can be understood as the mass and stiffness matrices of the anisotropic elliptic equation
\begin{equation}
- \nabla (\widehat{H}(x,y) \nabla u(x,y)) + L(y) u(x,y) = f(x,y), \qquad (x,y) \in \Omega,
\end{equation}
with $\widehat{H}$ as defined in \eqref{hat-H}. Hence, by following \cite{nie1985lumped}, the lumped counterparts $\widehat{M}$ of $\widetilde{M}$ and $\widehat{A}$ of $\widetilde{A}$ are defined by
\begin{align}
\label{x-normal-mass-lumped}
\begin{split}
\widehat{m}_{ij} &= \iint_\Omega I_h\left(\phi_i(x,y)\phi_j(x,y) L(y)\right)\dx\dy;\\
\widehat{a}_{ij} &= \iint_\Omega I_h\left(\left(\widehat{H}(x,y) \nabla\phi_i(x,y)\right) \cdot \nabla\phi_j(x,y)\right)\dx\dy,
\end{split}
\end{align}
for $i,j=1,\dots,(N+1)^2$, where $I_h$ is the element-wise, vector-valued interpolant operator, see \cite{frittelli2019preserving}. So, the tensorial decompositions of \eqref{x-normal-mass-lumped} are carried out as follows
\begin{equation}
\label{x-normal-decomposition-mass-lumped}
\begin{split}
\widehat{m}_{ij} &= \int_0^1 I_h(\psi_{i_x}(x)\psi_{j_x}(x))\dx \int_0^1 I_h(\psi_{i_y}(y)\psi_{j_y}(y)L(y))\dy = m^0_{i_xj_x} d^1_{i_yj_y} m^0_{i_yj_y};
\end{split}
\end{equation}
\begin{equation}
\label{x-normal-decomposition-stiffness-lumped}
\begin{split}
\widehat{a}_{ij} &= \int_0^1 \tdx{\psi_{ix}}(x)\tdx{\psi_{jx}}(x)\dx\int_0^1I_h\left(\psi_{i_y}(y)\psi_{j_y}(y)\frac{1}{L(y)}\right)\dy\\
&+ \int_0^1 I_h\left(\psi_{ix}(x)\psi_{jx}(x)\right)\dx\int_0^1I_h\left(\tdy{\psi_{i_y}}(y)\tdy{\psi_{j_y}}(y)L(y)\right)\dy\\
&+ \int_0^1 I_h\left(\tdx{\psi_{ix}}(x)\tdx{\psi_{jx}}(x)x^2\right)\dx\int_0^1I_h\left(\psi_{i_y}(y)\psi_{j_y}(y)\frac{L'^2(y)}{L(y)}\right)\dy\\
&- \int_0^1 I_h\left(\tdx{\psi_{ix}}(x)\psi_{jx}(x)x\right)\dx\int_0^1 I_h\left(\psi_{i_y}(y)\tdy{\psi_{j_y}}(y)L'(y)\right)\dy\\
&- \int_0^1 I_h\left(\psi_{ix}(x)\tdx{\psi_{jx}}(x)x\right)\dx\int_0^1 I_h\left( \tdy{\psi_{i_y}}(y)\psi_{j_y}(y)L'(y)\right)\dy\\
&= \frac{a_{i_xj_x} m_{i_yj_y}^{0}}{d^1_{i_y i_y}} + m^0_{i_xj_x} a_{i_yj_y}^{1} + \frac{a_{i_xj_x}^{2} m_{i_yj_y}^0 (d^2_{i_y i_y})^2}{d^1_{i_y i_y}} - c_{i_xj_x}d^3_{j_xj_x} c_{j_yi_y}d^2_{i_yi_y} - c_{j_xi_x}d^3_{i_x i_x} c_{i_yj_y}d^2_{j_y j_y}.
\end{split}
\end{equation}
Equations \eqref{x-normal-decomposition-mass-lumped}-\eqref{x-normal-decomposition-stiffness-lumped} translate to the following matrix identities:
\begin{align}
\label{kronecker_mass_distorted_lumped}
\begin{split}
\widehat{M} &= M_0\otimes M_0 D_1;\\
\widehat{A} &= A\otimes D_1^{-1} M_0+ M_0\otimes A_1 + A_2 \otimes D_2^{2} D_1^{-1} M_0 - C D_3\otimes D_2 C^T - D_3 C^T \otimes C D_2.
\end{split}
\end{align}
As mentioned earlier, the matrix relations \eqref{kronecker_mass_distorted_lumped} are simpler than \eqref{kronecker_mass_distorted}. This is because (i) the matrices $D_1, D_2, D_3, M_0$ are now diagonal and (ii) the only non-diagonal matrices involved, i.e. $A,A_1,A_2,C$ are now tridiagonal.

\begin{remark}[Homogeneous Dirichlet boundary conditions]
\label{rmk:boundary_cond}
In the presence of homogeneous Dirichlet boundary conditions, it is well-known that the boundary basis functions must be eliminated from the basis, see \cite{hughes2012finite}, hence all the matrices involved in \eqref{kronecker_mass_distorted_lumped} must be trimmed by eliminating all boundary entries. Hence, the dimension of such matrices drops from $(N+1)\times (N+1)$ to $(N-1)\times (N-1)$ and the following additional properties hold true.
\begin{itemize}
\item $M_0$ is a multiple of the identity: for $m_0 := \frac{1}{N}$ it holds that
\begin{equation}
\label{dirichlet_lumped_constant}
M_0 = m_0 I.
\end{equation}
\item The stiffness matrix $A$ is Toeplitz.
\item Thanks to the symmetries of the $\psi_i$'s, $C$ is Toeplitz and skew-symmetric, i.e.
\begin{equation}
\label{dirichlet_convection_skew_symmetric}
C^T = -C.
\end{equation}
\end{itemize}
\end{remark}

\noindent
Because the dimension of the matrices depends on the kind of boundary conditions, see Remark \ref{rmk:boundary_cond}, we set
\begin{equation}
\label{variable_dimension}
q :=
\begin{cases}
N+1 \qquad \text{for Neumann boundary conditions}\\
N-1 \qquad \text{for Dirichlet boundary conditions}.
\end{cases}
\end{equation}
In the presence of homogeneous Dirichlet boundary conditions, thanks to \eqref{dirichlet_lumped_constant} and \eqref{dirichlet_convection_skew_symmetric}, relations \eqref{kronecker_mass_distorted_lumped} become
\begin{align}
\label{kronecker_mass_distorted_lumped_dirichlet}
\begin{split}
&\widehat{M} = m_0^2 I\otimes  D_1;\\
&\widehat{A} = m_0(A\otimes D_1^{-1} + I \otimes A_1 + A_2 \otimes D_2^{2} D_1^{-1}) + C D_3\otimes D_2 C + D_3 C \otimes C D_2.
\end{split}
\end{align}

\begin{remark}[Recap on lumped $\mathbb{P}_1$ matrix properties]
We recap here the properties of the matrices appearing in \eqref{kronecker_mass_distorted_lumped}and \eqref{kronecker_mass_distorted_lumped_dirichlet}:
\begin{itemize}
\item $M_0$ and $D_1$ are diagonal and positive definite. In the Dirichlet case $M_0$ is multiple of the identity;
\item $D_2$ is diagonal and it is non-singular only when the curve $x = L(y)$ is strictly monotone;
\item $D_3$ is diagonal and is singular when $N$ is odd, or tends to being singular when $N$ is even and approaches infinity;
\item $A, A_1$ are tridiagonal, moreover they are positive definite in the Dirichlet case and semidefinite in the Neumann case (one null eigenvalue). $A$ is symmetric and in the Dirichlet case it is Toeplitz;
\item $A_2$ is tridiagonal non-symmetric, it is singular for $N$ odd, or tends to being singular for $N$ even approaching infinity;
\item $C=\texttt{tridiag}([-1, 0, 1])$, except for boundary entries in the Neumann case (it is the matrix of centered first derivatives). According to $N$ being even or odd, and depending on the boundary conditions, $C$ is singular or tends to being singular as $N \rightarrow + \infty$.
\end{itemize}
\end{remark}

\section{Stationary PDEs and the Sylvester equation}
\label{sec:poisson_square}
In this section we consider the Poisson equation on the unit square $\Omega = [0,1]^2$:
\begin{equation}
\label{poisson_equation_square}
-\Delta u(\boldx) + \gamma u(\boldx) = f(\boldx), \qquad \boldx\in\Omega,
\end{equation}
where $\gamma \geq 0$ for the case of Dirichlet boundary conditions and $\gamma > 0$ for the case of Neumann boundary conditions. Since we encompass general boundary condition, we use the notation in \eqref{variable_dimension} for the dimension of matrices. The general FEM discretisation of problem \eqref{poisson_equation_square} in vector form is then
\begin{equation}
\label{poisson_square_FEM}
\widetilde{A}\boldxi + \gamma \widetilde{M}\boldxi = \widetilde{M}\boldf,
\end{equation}
with $\boldxi$ being the nodal vector of the numerical solution and $\boldf$ being the corresponding nodal vector of $f$. By using \eqref{kronecker_mass}, the linear system \eqref{poisson_square_FEM} becomes the following generalised Sylvester matrix equation:
\begin{equation}
\label{laplace_dirichlet_sylvester}
\left(A+\frac{\gamma}{2}M\right)UM + MU\left(A+\frac{\gamma}{2}M\right) = MFM,
\end{equation}
where $U$ and $F$ are such that $\texttt{vec}(U) = \boldxi$ and $\texttt{vec}(F) = \boldf$. Since the mass matrix $M$ is positive-definite, we can pre- and post-multiply both sides of \eqref{laplace_dirichlet_sylvester} by $M^{-1}$:
\begin{equation}
\label{neumann_sylvester_0}
\left(M^{-1}A+\frac{\gamma}{2}I\right)U + U\left(AM^{-1}+\frac{\gamma}{2}I\right) = F.
\end{equation}
Even if $M^{-1}A$ and $AM^{-1}$ are not symmetric, they are diagonalizable nonetheless because $A,M$ are both symmetric and $M$ is positive definite. Hence, also $Z_1 := M^{-1}A + \frac{\gamma}{2} I$ and $Z_2 := AM^{-1} + \frac{\gamma}{2} I$ are diagonalizable. Then, \eqref{neumann_sylvester_1} becomes
\begin{equation}
\label{neumann_sylvester_1}
Z_1U + UZ_2 = F,
\end{equation}
and we can diagonalise $Z_1$, $Z_2$ as follows
\begin{align}
\label{diagonalisation_1}
Z_1 = X^{(1)}\Lambda^{(1)} {X^{(1)}}^{-1}; \qquad Z_2 = X^{(2)}\Lambda^{(2)} {X^{(2)}}^{-1},
\end{align}
with $\Lambda^{(k)} \in\mathbb{R}^{q\times q}$, $k=1,2$, being the diagonal matrices containing the eigenvalues $\lambda_{i}^{(k)}$, $i=1,\dots,q$ of $Z_1$ and $Z_2$, respectively. Except special cases, the $\Lambda^{(k)}$'s and $X^{(k)}$'s must be computed numerically. Now, by setting $\widehat{U} := {X^{(1)}}^{-1} U X^{(2)}$ and $\widehat{F} := {X^{(1)}}^{-1} F X^{(2)}$, \eqref{neumann_sylvester_1} becomes
\begin{equation}
\label{neumann_sylvester_2}
\Lambda^{(1)} \widehat{U} + \widehat{U} \Lambda^{(2)} = \widehat{F},
\end{equation}
which can be solved as follows. Let $L\in\mathbb{R}^{q\times q}$ be the matrix defined by
\begin{equation}
\ell_{ij} = \frac{1}{\lambda_{i} + \lambda_{j}}, \qquad i,j=1,\dots, q.
\end{equation}
As shown in \cite{dautilia2020matrix}, the solution to \eqref{neumann_sylvester_2} in the spectral space can be expressed as:
\begin{equation}
\widehat{U} = L \circ \widehat{F},
\end{equation}
with $\circ$ denoting the Hadamard product. The original variable $U$ is thus given by
\begin{equation}
\label{reduced_formula_for_u}
U = X^{(1)} (L \circ \widehat{F}) {X^{(2)}}^{-1}.
\end{equation}
We call this technique \emph{FEM reduced method}.
\begin{remark}[Memory performance of the reduced approach]
\label{rmk:memory_performance_reduced}
Because $U$, $F$, $L$, $X^{(1)}$, $X^{(2)}$ are the only full matrices involved in the computation, the memory occupation of the proposed approach is $5N^2 + O(kN)$ floating point numbers for any $k$. On the other hand, since the large matrices $\widetilde{A}$ and $\widetilde{M}$ are $(2k+1)^2$-diagonal, where $k$ is the polynomial order of the method, and $\boldf$ and $\boldxi$ are both full vectors, the vector formulation \eqref{poisson_square_FEM} has a memory occupation of $((2k+1)^2+2)N^2 + O(kN)$ floating point numbers.
\end{remark}

\subsection{Special case: Lagrangian $\mathbb{P}_1$ elements with Dirichlet boundary conditions}
In the special case of Lagrangian $\mathbb{P}_1$ elements with Dirichlet boundary conditions, we are able to solve the matrix equation \eqref{neumann_sylvester_1} in closed form without computing spectral decompositions numerically. In fact, since $M$ and $A$ are both symmetric, positive definite, tridiagonal Toeplitz matrices, there exist $\alpha,\beta>0$ such that $M = \alpha A + \beta I$. In this specific case, we have
\begin{equation}
\label{dirichlet_linear_combination}
\alpha = -\frac{1}{6N^2}, \qquad \beta = \frac{1}{N}.
\end{equation}
This implies that (i) $M^{-1}$ and $A$ commute, (ii) $M^{-1}$ and $A$ share the same eigenvectors and (iii) $M^{-1}A = AM^{-1}$ is diagonalizable and shares the same eigenvectors of $M^{-1}$ and $A$. We can thus write
\begin{align}
&Z_1 = Z_2 = X\Lambda X^{-1},
\end{align}
with $\Lambda \in\mathbb{R}^{(N-1)\times (N-1)}$ being the diagonal matrix containing the eigenvalues $\lambda_i$, $i=1,\dots,N-1$ of $Z_1 = Z_2$. Now, by setting $\widehat{U} := X^{-1} U X$ and $\widehat{F} := X^{-1} F X$, \eqref{neumann_sylvester_1} becomes
\begin{equation}
\label{dirichlet_sylvester_2}
\Lambda \widehat{U} + \widehat{U} \Lambda = \widehat{F},
\end{equation}
in the spectral space. We are left to show that $X$ and $\Lambda$ are known in closed form. In fact, as shown in \cite{chung2000discrete}, the eigenvalues and eigenvectors of $A$ are given by
\begin{align}
\label{dirichlet_eigenvalues}
\lambda^A_{i} = \frac{2}{N}\left(1-\cos\frac{i\pi}{N} \right), \qquad i=1,\dots,N-1;\\
\label{dirichlet_eigenvectors}
(\boldv^A_{i})_j = \sin\frac{ij\pi}{N}, \qquad i,j = 1,\dots,N-1.
\end{align}
Hence, by using \eqref{dirichlet_linear_combination} the entries of $\Lambda$ are given by
\begin{equation}
\label{dirichlet_combined_eigenvalues}
\lambda_i = \frac{(12N^2-\gamma)\lambda^A_i + 6N\gamma}{12N - 2\lambda_i^A}, \quad i=1,\dots, N-1,
\end{equation}
while the entries of $X$ (common basis of eigenvectors of $A$, $Z_1$ and $Z_2$) are given in \eqref{dirichlet_eigenvectors}.  Hence,  $U$ is given by \eqref{reduced_formula_for_u} with $\Lambda^{(1)} = \Lambda^{(2)} = \Lambda$ defined by \eqref{dirichlet_combined_eigenvalues} and $X^{(1)} = X^{(2)} = X$ defined by \eqref{dirichlet_eigenvectors}.

\subsection{Numerical example}
\label{sec:example_poisson_reduced_pk}
\begin{figure}[t!]
\begin{center}
\hspace*{-20mm}
\includegraphics[scale=0.45]{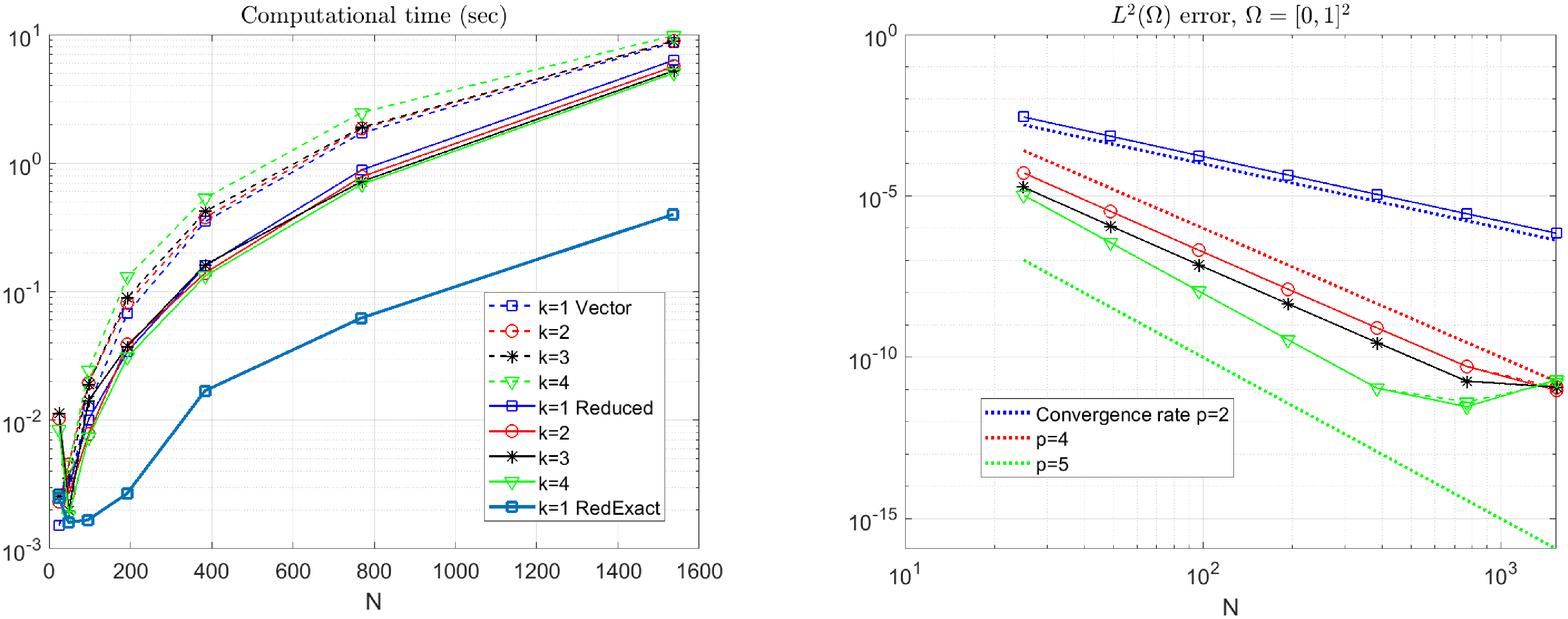}
\end{center}
\caption{Poisson equation \eqref{example_laplace_square_dirichlet} on the square $\Omega=[0,1]^2$. FEM of order  $k=1,2,3,4$: comparison between the  Kronecker (vector) approach \eqref{poisson_square_FEM} (dashed lines), the MO-FEM solved by the \emph{reduced} approach \eqref{neumann_sylvester_2} (continuous lines) and, only for $k=1$, the reduced approach in closed form \eqref{dirichlet_sylvester_2} (thick continuous line). Left plot: computational times. Right plot: convergence behavior. Dotted lines indicate slopes 2,4,5.}
\label{fig:square_dirichlet_convergence_pk}
\end{figure}

We consider the following Poisson equation with Dirichlet boundary conditions on the square $\Omega = [0,1]^2$:
\begin{equation}
\label{example_laplace_square_dirichlet}
\begin{cases}
-&\Delta u(x,y) = 8\pi^2\sin(2\pi x)\sin(2\pi y), \qquad (x,y) \in \Omega;\\
&u(x,y) = 0, \qquad (x,y) \in  \partial \Omega,
\end{cases}
\end{equation}
whose exact solution is $u(x,y) = \sin(2\pi x)\sin(2\pi y)$. We apply Lagrangian $\mathbb{P}_k$ elements, $k=1,\dots,4$, and we compare the performances of the Kronecker (vector) formulation \eqref{poisson_square_FEM}, solved through the direct solver \texttt{mldivide} of MATLAB (known also as "backslash" $\backslash$ ), and the \emph{reduced} approach \eqref{neumann_sylvester_2}. For $k=1$, we further compare the aforementioned methods with the \emph{reduced} method in closed form \eqref{dirichlet_sylvester_2}. Here we consider a sequence of $7$ meshes $\Gamma_i$ with $N_i = 24\cdot 2^i$ for all $i=0,\dots,6$. Such $N_i$'s are compatible with $\mathbb{P}_k$ elements for all $k=1,\dots,4$.  The numerical results are shown in Fig. \ref{fig:square_dirichlet_convergence_pk}.
On the finest mesh ($N=1536$), the reduced approach \eqref{neumann_sylvester_2} for $k=1$ is approximately $1.38$ times quicker than the direct (vector) method, $1.58$ times quicker for $k=2$, $1.70$ times quicker for $k=3$, and $1.96$ times quicker for $k=4$. The \emph{reduced} approach in closed form \eqref{dirichlet_sylvester_2} (only $k=1$) is $21.73$ times quicker than the vector formulation. Furthermore, on equal meshes, the methods produce the same solutions up to rounding errors (they are equivalent to each other) and exhibit optimal convergence in space, with the case $k=2$ being superconvergent (fourth order), as we can see in Fig. \ref{fig:square_dirichlet_convergence_pk}, right plot.\\
This and all the following experiments are carried out in MATLAB R2019a on a HP Z230 Tower Workstation with Intel Core i7-440 CPU and 16GB RAM. The timings were just taken once and not averaged
with several measurements.

\section{Stationary PDEs on $x$-normal domains and multiterm Sylvester equations}
\label{sec:implementation-cap-shaped-domain-dirichlet-lumped}

We consider the following stationary PDE problem on an $x$-normal domain $\Omega^L$:
\begin{equation}
\label{laplace_equation_xnormal}
- \Delta u(\boldx) + \gamma u(\boldx) = f(\boldx), \qquad \boldx \in \Omega^L,
\end{equation}
with $\gamma \geq 0$ in the presence of zero Dirichlet boundary conditions and $\gamma>0$ in the presence of zero Neumann boundary conditions. The FEM discretisation of problem \eqref{laplace_equation_xnormal} in vector form is then
\begin{equation}
\label{laplace_xnormal_FEM}
\widetilde{A}\boldxi + \gamma\widetilde{M}\boldxi = \widetilde{M}\boldf,
\end{equation}
with $\widetilde{A}$ and $\widetilde{M}$ as defined in \eqref{x-normal-stiffness-mass}. The lumped counterpart of \eqref{laplace_xnormal_FEM} is
\begin{equation}
\label{laplace_xnormal_FEM_lumped}
\widehat{A}\boldxi + \gamma\widehat{M}\boldxi = \widehat{M}\boldf,
\end{equation}
with $\widehat{M}$ and $\widehat{A}$ in \eqref{x-normal-mass-lumped}. Here, $\boldxi \in\mathbb{R}^{q^2}$ is the nodal vector of the numerical solution and $\boldf\in\mathbb{R}^{q^2}$ being the corresponding nodal vector of $f$. By using \eqref{kronecker_mass_distorted}, the linear system \eqref{laplace_xnormal_FEM} becomes the following multiterm Sylvester matrix equation:
\begin{equation}
\label{laplace_cap-shaped-domain_sylvester}
M_1UA + (B_1 + \gamma M_3)UM + M_2UB_2 - C_2UC_1 - C_2^TUC_1^T = M_3FM.
\end{equation}
Similarly, by using \eqref{kronecker_mass_distorted_lumped_dirichlet} the ``lumped'' linear system \eqref{laplace_xnormal_FEM_lumped} translates to:
\begin{equation}
\label{laplace_cap-shaped-domain_sylvester_lumped}
\begin{split}
M_0D_1^{-1}UA + (A_1+\gamma M_0D_1)UM_0 + M_0D_1^{-1}D_2^2UA_2 - CD_2UCD_3 &- D_2C^TUD_3C^T\\
&= M_0D_1FM_0.
\end{split}
\end{equation}
Thanks to \eqref{dirichlet_convection_skew_symmetric}, in the case of zero Dirichlet boundary conditions, \eqref{laplace_cap-shaped-domain_sylvester_lumped} can be simplified as:
\begin{equation}
\label{laplace_cap-shaped-domain_sylvester_dirichlet_lumped}
\begin{split}
m_0(D_1^{-1}UA + (A_1 +\gamma D_1)U + D_2^2D_1^{-1}UA_2) - (D_2CUD_3C + CD_2UCD_3) = m_0^2D_1F.
\end{split}
\end{equation}
We now propose an iterative strategy that can  be applied for the solution of the multiterm matrix equations \eqref{laplace_cap-shaped-domain_sylvester}, \eqref{laplace_cap-shaped-domain_sylvester_lumped} and \eqref{laplace_cap-shaped-domain_sylvester_dirichlet_lumped}.  For the systems \eqref{laplace_xnormal_FEM} and \eqref{laplace_xnormal_FEM_lumped}, the well known Preconditioned Conjugate Gradient method (PCG) \cite{golub2013matrix} would be a suitable choice since the matrices $\widetilde{A} +\gamma \widetilde{M}$ and $\widehat{A} + \gamma \widehat{M}$ are symmetric and positive definite. Here, we propose its matrix oriented version to solve the corresponding Sylvester multiterm equations \eqref{laplace_cap-shaped-domain_sylvester} and \eqref{laplace_cap-shaped-domain_sylvester_dirichlet_lumped}. To this end, we define the following matrix operators:
\begin{align}
&\widetilde{\mathcal{L}}(U) := M_1UA + (B_1 + \gamma M_3)UM + M_2UB_2 - C_2UC_1 - C_2^TUC_1^T;\\
&\widetilde{\mathcal{R}}(U) := M_3FM;\\
&\widehat{\mathcal{L}}(U) := M_0D_1^{-1}UA + (A_1+\gamma M_0D_1)UM_0 + M_0D_1^{-1}D_2^2UA_2 - CD_2UCD_3 - D_2C^TUD_3C^T;\\
&\widehat{\mathcal{R}}(U) := M_0D_1FM_0;
\end{align}
for all $U\in\mathbb{R}^{q\times q}$. We also consider suitable preconditioning operators $\widetilde{\mathcal{P}}, \widehat{\mathcal{P}}: \mathbb{R}^{q\times q} \rightarrow \mathbb{R}^{q\times q}$ for the FEM and lumped FEM, respectively, whose choice will be discussed later. With these settings, the matrix-oriented formulation of the PCG for systems \eqref{laplace_xnormal_FEM} and \eqref{laplace_xnormal_FEM_lumped}, that we define as MO-PCG, is given by
\begin{align}
\label{laplace_xnormal_pcg}
&\begin{cases}
R^{(0)} = \mathcal{R}(U^{(0)}) - \mathcal{L}(U^{(0)});\\
Z^{(0)} = \mathcal{P}^{-1}(R^{(0)});\\
Q^{(0)} = Z^{(0)};
\end{cases}\\
\label{laplace_xnormal_pcg_iteration}
&\begin{cases}
\alpha^{(s)} = \dfrac{\tsum(\tsum(Q^{(s)}\circ R^{(s)}))}{\tsum(\tsum(\mathcal{L}(Q^{(s)})\circ R^{(s)}))};\\
U^{(s+1)} = U^{(s)} + \alpha^{(s)}Q^{(s)};\\
R^{(s+1)} = R^{(s)} - \alpha^{(s)}\mathcal{L}(Q^{(s)});\\
Z^{(s+1)} = \mathcal{P}^{-1}(R^{(s+1)});\\
\beta^{(s)} = \dfrac{\tsum(\tsum(\mathcal{L}(Q^{(s)})\circ Z^{(s+1)}))}{\tsum(\tsum(\mathcal{L}(Q^{(s)})\circ R^{(s)}))};\\
Q^{(s+1)} = Q^{(s)} - \beta^{(s)}Q^{(s)},
\end{cases}
\qquad  s \geq 0,
\end{align}
where the $\ \widetilde{}\ $ and $\ \widehat{}\ $ are omitted for ease of presentation. For the initial guess $U^{(0)}$ we can choose for instance $U^{(0)} = \boldzero$.\\
For illustrative purposes, we consider the following stopping criterion. If $err_v$ is the absolute error obtained by the direct method solving the linear systems (vector formulation) \eqref{laplace_xnormal_FEM} or \eqref{laplace_xnormal_FEM_lumped},  we stop the MO-PCG iterations when
the increment fulfils
$$\|U^{(s+1)} - U^{(s)}\|_F \leq 0.05 \| err_v \|_F $$ where $\|\cdot\|_F$ denotes the Frobenius norm.

\subsubsection*{Memory performance of the PCG approach}
Because $U^{(s)}$, $R^{(s)}$, $Z^{(s)}$, $Q^{(s)}$, $F$ are the only full matrices involved in the computation, the memory occupation of the proposed approach is $5N^2 + O(kN)$, with $N$ being the grid size and $k$ is the FEM polynomial order On the other hand, the vector - Kronecker formulation \eqref{poisson_square_FEM} has a memory occupation of $((2k+1)^2+2)N^2 + O(kN)$  for matrix storage, where $k$ is the polynomial order of the method, as discussed in Remark \ref{rmk:memory_performance_reduced}.

\subsubsection*{Matrix PCG vs. classical PCG}
The large linear system \eqref{laplace_xnormal_FEM} or \eqref{laplace_xnormal_FEM_lumped} could be solved via PCG in classical vector formulation. However, in the experiments (carried out in MATLAB R2019a on a HP Z230 Tower Workstation with Intel Core i7-4770 CPU and 16 GB RAM), we find that the matrix-oriented PCG is significantly faster than its classical vector counterpart, even if the methods are equivalent. Some details are provided in the test of Section \ref{sec:experiment_elliptic_cap-shaped-domain_curved}. 

\subsubsection*{Choice of the preconditioners}
The choice of fast and efficient matrix-oriented preconditioners for problems \eqref{laplace_cap-shaped-domain_sylvester}-\eqref{laplace_cap-shaped-domain_sylvester_lumped} is an open problem. Here we consider the case of Dirichlet boundary conditions and we derive preconditioners that experimentally prove more efficient than the identity operator, i.e. no preconditioner.
To this end, we first construct suitable preconditioners for the large linear systems \eqref{laplace_xnormal_FEM}-\eqref{laplace_xnormal_FEM_lumped}, which we will use to derive matrix-oriented preconditioners for the multiterm Sylvester problems \eqref{laplace_cap-shaped-domain_sylvester}-\eqref{laplace_cap-shaped-domain_sylvester_dirichlet_lumped}.
For the systems \eqref{laplace_xnormal_FEM}-\eqref{laplace_xnormal_FEM_lumped} in matrix form, we could consider the \emph{ideal} preconditioners
\begin{align}
\label{ideal_preconditioner}
&\widetilde{P} := A \otimes M_1 + M \otimes B_1 + B_2 \otimes M_2; \qquad \widehat{P} := A\otimes D_1^{-1}M_0 + M_0 \otimes A_1 + A_2 \otimes D_2^{2} D_1^{-1}M_0,
\end{align}
which contain the discrete operators for the second-order derivative terms of the stiffness matrices $\widetilde{A}$ and $\widehat{A}$, respectively. Nevertheless, in matrix form the preconditioning operations $P^{-1}\texttt{vec}(U)$, with $P$ equal $\widetilde{P}$ or $\widehat{P}$, imply again the solution of another multiterm Sylvester equation at each iteration of the PCG due to the presence of more Kronecker products in \eqref{ideal_preconditioner}. In order to avoid this, we approximate in a spectral sense the preconditioners $\widetilde{P}$ and $\widehat{P}$ in \eqref{ideal_preconditioner}, respectively, such that a single Kronecker product is present. After several experiments, the choice $\widetilde{P}_1 := (A+B_2) \otimes B_1$ and $\widehat{P}_1 := (A+A_2) \otimes A_1$ proved the best among the tested ones. In matrix form, these preconditioners translate to the following operators:
\begin{align}
\label{vector_preconditioner}
\widetilde{\mathcal{P}}(U) := B_1U(A+B_2); \qquad \widehat{\mathcal{P}}(U) := A_1U(A+A_2).
\end{align}
that satisfy  $\texttt{vec}(\widetilde{\mathcal{P}}(U)) = \widetilde{P}_1\texttt{vec}(U)$ and $\texttt{vec}(\widehat{\mathcal{P}}(U)) = \widehat{P}_1\texttt{vec}(U)$ for all $U\in\mathbb{R}^{q\times q}$. The respective inverse operators are given by
$$\widetilde{\mathcal{P}}^{-1}(U) = B_1^{-1}U(A+B_2)^{-1};\qquad \widehat{\mathcal{P}}^{-1}(U) = A_1^{-1}U(A+A_2)^{-1}.$$
In the special case of square domains, the choices \eqref{vector_preconditioner} reduce to
\begin{align}
\label{precSquare}
\widetilde{\mathcal{P}}(U) = \widehat{\mathcal{P}}(U) = AUA,
\end{align}
which will prove particularly efficient in the following experiments. 

\begin{figure}[t!]
\begin{center}
\hspace*{-10mm}
\includegraphics[scale=0.45]{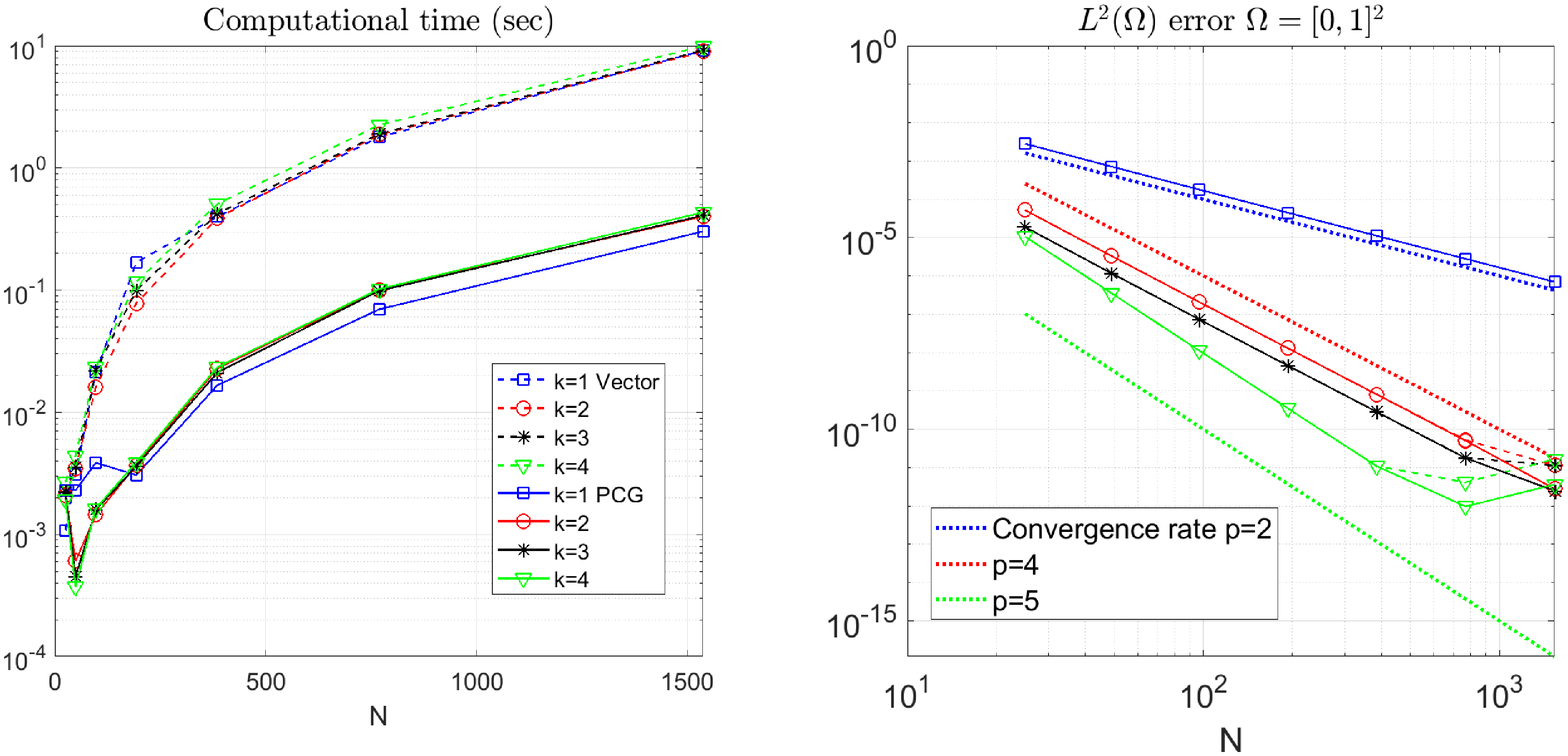}
\end{center}
\caption{Poisson equation \eqref{example_laplace_square_dirichlet} on the square $\Omega=[0,1]^2$.  FEM of order  $k=1,2,3,4$: comparison between the  "vector"  formulation \eqref{poisson_square_FEM} (dashed lines) solved by the direct method and the MO-FEM solved by the MO-PCG  (continuous lines). Left plot: computational times. Right plot: Convergence behaviour of the two approaches for all $k$. Dotted lines indicate slopes 2,4 and 5. For $k=1,3,4$, the convergence is optimal of order $k+1$, for $k=2$, we observe superconvergence of order 4. PCG is cheaper in execution time and for all $k$ and for all $N$ stops after two iterations, see also Table 1.}
\label{fig:elliptic_square_PCG}
\end{figure}

\subsection{Example 1: Poisson equation on the square, $\mathbb{P}_k$ elements}
We consider the Poisson equation \eqref{example_laplace_square_dirichlet} on the square $\Omega = [0,1]^2$. We consider Lagrangian $\mathbb{P}_k$ finite elements for $k=1,2,3,4$, and we compare the classical vector approach that solves the Kronecker form \eqref{poisson_square_FEM} through the direct solver \texttt{mldivide} of MATLAB, and the MO-PCG approach \eqref{laplace_xnormal_pcg}-\eqref{laplace_xnormal_pcg_iteration} using \eqref{precSquare} as preconditioner. Observe that, since the domain is the unit square, the multiterm Sylvester formulation \eqref{laplace_cap-shaped-domain_sylvester} reduces to the classical Sylvester equation \eqref{laplace_dirichlet_sylvester} and then in this case the MO-PCG can be regarded also as an alternative to the \emph{reduced} approach presented in Section \ref{sec:poisson_square}. We therefore consider the same sequence of seven meshes considered in Section \ref{sec:example_poisson_reduced_pk} with $N=24\cdot 2^i$, $i=0,\dots,6$. On equal meshes, the MO-PCG with preconditioning in \eqref{precSquare} is quicker than the direct method for the vector approach and the gap increases with $N$, as we can see in Fig. \ref{fig:elliptic_square_PCG}, left plot. 
The detailed time comparisons on the finest mesh are shown in Table \ref{tab:square_pcg_times}. For all $k$, the methods exhibit optimal convergence ($(k+1)$-th order), with the case $k=2$ being superconvergent (fourth order) as we can see in Fig. \ref{fig:elliptic_square_PCG}, right plot. On the finest mesh ($N=1536, k=4$) errors near the machine precision affect the convergence order, with MO-PCG being more accurate. For all $k$ and for all $N$, PCG terminates with two iterations. In conclusion, the MO-PCG method outperforms also the \emph{reduced matrix} approach \eqref{neumann_sylvester_2} when the eigenvalue decompositions \eqref{diagonalisation_1} are computed numerically.

\begin{table}[t!]
\caption{Poisson equation \eqref{example_laplace_square_dirichlet} on the square $\Omega=[0,1]^2$: computational times on the finest mesh ($N = 1536$) for all $k=1,2,3,4$ (see Fig.4, left plot) and the respective time ratios between the direct method for the vector form \eqref{poisson_square_FEM} and the MO-PCG method \eqref{laplace_xnormal_pcg}-\eqref{laplace_xnormal_pcg_iteration} with preconditioning in \eqref{precSquare}.}
\begin{center}
\begin{tabular}{ l c c c }
$k$ & \thead{Time (s)\\   Direct method (vector form)} & \thead{Time (s)\\ Matrix PCG} & \thead{Time ratio \\ (Direct/PCG)}\\ 
 \hline
 $1$  & 9.0273 & 0.2992 & 30.1726\\  
 $2$  & 9.4791 & 0.4255 & 22.2760\\
 $3$  & 9.3019 & 0.4354 & 21.3665\\
 $4$  & 10.1502 & 0.4438 & 22.8726\\
\end{tabular}
\end{center}
\label{tab:square_pcg_times}
\end{table}

\subsection{Example 2: Poisson equation on curved domain}\label{sec:experiment_elliptic_cap-shaped-domain_curved}

\begin{figure}[t!]
\begin{center}
\includegraphics[angle=0,scale=0.6]{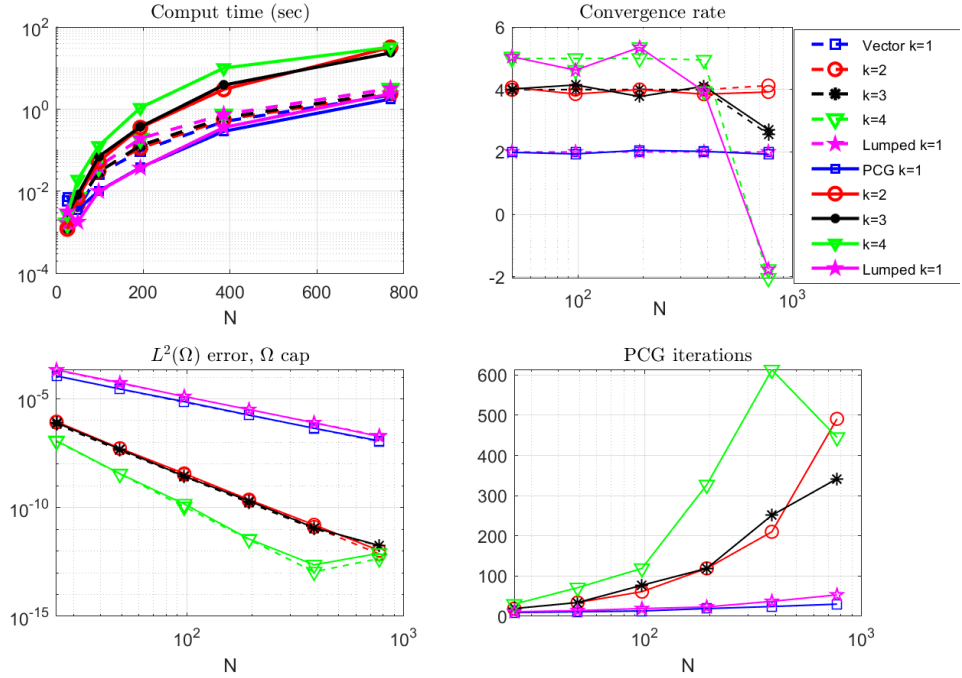}
\end{center}
\caption{Poisson equation \eqref{example_laplace_cap-shaped-domain_dirichlet} on the cap-shaped domain $\Omega^S$ defined in \eqref{cap-shaped-domain}. Continuous lines indicate MO-PCG, while dashed lines indicate the vector (Kronecker) formulations solved by the direct method. In terms of computational time, MO-PCG is competitive only for $k=1$, {without lumping} (upper left panel). Both approaches exhibit optimal convergence for all FEM order $k$, with the case $k=2$ being superconvergent, round-off error becomes dominant on the finest mesh for $k=4$ (upper-right and lower-left panels). Number of iterations required by the MO-PCG increases with $k$ and $N$, lower-right panel.}
\label{fig:poisson_pk}
\end{figure}

Consider the \emph{cap-shaped} symmetric $x$-normal domain
\begin{equation}
\label{cap-shaped-domain}
\Omega^S := \left\{(x,y)\in\mathbb{R}^2 \middle | 0 \leq y \leq 1, \ |x| \leq 1-\frac{y^2}{2}\right\},
\end{equation}
(shown in the next Fig. \ref{fig:spot_laby}). On $\Omega^S$ we consider the Poisson equation with zero Dirichlet boundary conditions:
\begin{equation}
\label{example_laplace_cap-shaped-domain_dirichlet}
\begin{cases}
-&\Delta u(x,y) = f(x,y), \qquad (x,y) \in \Omega^S;\\
&u(x,y) = 0, \qquad (x,y) \in  \partial \Omega^S,
\end{cases}
\end{equation}
were $f(x,y)$ is chosen in such a way that the exact solution is 
\begin{equation*}
u(x,y) = y(y - 1)\left(- \frac{y^2}{2} + x + 1\right)\left(\frac{y^2}{2} + x - 1\right), \qquad (x,y) \in \Omega^S,
\end{equation*}
we omit the cumbersome expression of such $f(x,y)$. 
We consider the $\mathbb{P}_k$ elements \eqref{laplace_xnormal_FEM}, $k=1,\dots,4$, and the lumped $\mathbb{P}_1$ elements in vector form \eqref{laplace_xnormal_FEM_lumped} and we solve again these classical linear systems via the MATLAB direct solver \texttt{mldivide}. Hence, we compare the perfomance in terms of computational time with the matrix-oriented PCG method \eqref{laplace_xnormal_pcg}-\eqref{laplace_xnormal_pcg_iteration} with preconditioners as in \eqref{vector_preconditioner}, respectively, that solve the multiterm Sylvester equations arising by the MO-FEM.

By applying all methods on a sequence of $6$ meshes $\Omega_i^S$, $i=0,\dots,5$ with $N_i = 24\cdot 2^i$ for all $i=0,\dots,5$ we find optimal quadratic convergence in $L^2(\Omega^S)$ and almost equal errors, as reported in Fig. \ref{fig:poisson_pk}, upper right plot. By comparing the computational times, as shown in Fig. \ref{fig:poisson_pk}, upper left plot, the MO-PCG is competitive only for $k=1$ and moderate meshsizes $N < 768$. Note that the no lumped FEM is slightly better than the lumped version. On finer meshes,  say for $N \geq 800$, the number of iterations required by MO-PCG increases dramatically for $k \geq 2$, so affecting its global performance.
We guess that different preconditioners could improve this behaviour, but this study is outside the scope of the present work.
On the other hand, we find in this particular experiment that, the PCG method in its classical vector form is much slower the our MO-PCG, as discussed before. In fact, by applying the built-in function \emph{pcg} of Matlab with the same preconditioner, we find that MO-PCG is much quicker and the gap increases with $N$. For example, for $k=1$ without lumping, the speedup factor is:  $1.48$ for $N=24$ and $30.82$ for $N = 768$.  As a final remark, it is worth noting that in any case, as $N$ and $k$ increase, MO-PCG becomes more competitive than the Direct (Kronecker) solver in terms of memory occupation, as explained before in more detail.

In the next section we will show that the MO-PCG approach will prove far more convenient in the case of time-dependent PDEs.

\section{Time-dependent PDEs}
\label{sec:time_dependent_pdes}
On a general $x$-normal domain $\Omega^L$ we consider the following semilinear heat equation
\begin{equation}
\label{semilinear_heat_equation}
u_t -d_u\Delta u = f(u,\boldx,t), \qquad \boldx \in \Omega^L, \ t \in [0,T],
\end{equation}
with $d_u>0$, endowed with either homogeneous Dirichlet boundary conditions $u_{|\partial \Omega^S} = 0$ or homogeneous Neumann boundary conditions $(\nabla u \cdot \boldn)_{|\partial \Omega^S} = 0$. The treatment of RDSs of the form \eqref{general_RDS} is completely analogous and we omit the details. The $\mathbb{P}_k$ and lumped $\mathbb{P}_1$ FEM spatial discretisations of \eqref{semilinear_heat_equation} in vector form are as follows
\begin{align}
\label{semilinear_heat_equation_semidiscrete_pk}
\widetilde{M}\boldxi_t + d_u\widetilde{A}\boldxi = \widetilde{M}f(\boldxi),\\
\label{semilinear_heat_equation_semidiscrete}
\widehat{M}\boldxi_t + d_u\widehat{A}\boldxi = \widehat{M}f(\boldxi),
\end{align}
respectively, where $\boldxi(t)\in\mathbb{R}^{q^2}$,  with $q$ as defined in \eqref{variable_dimension}, is the time-dependent nodal vector of the spatially discrete solution. A full discretisation can be obtained by applying the IMEX Euler timestepping scheme with timestep $\tau>0$ to \eqref{semilinear_heat_equation_semidiscrete_pk} and \eqref{semilinear_heat_equation_semidiscrete}, which yields
\begin{align}
\label{semilinear_heat_equation_fully_discrete_pk}
(I + d_u\tau\widetilde{M}^{-1}\widetilde{A})\boldxi^{(k+1)} = \boldxi^{(k)} + \tau \boldf^{(k)}, \qquad k = 0, \dots, N_T-1,\\
\label{semilinear_heat_equation_fully_discrete}
(I + d_u\tau\widehat{M}^{-1}\widehat{A})\boldxi^{(k+1)} = \boldxi^{(k)} + \tau \boldf^{(k)}, \qquad k = 0, \dots, N_T-1,
\end{align} 
respectively, where $N_T := \left\lceil\frac{T}{\tau}\right\rceil$, $\boldxi^{(k)}$ is the nodal vector of the fully discrete solution at time $t_k := k\tau$ and $\boldf^{(k)} := f(\boldxi^{(k)})$, see \cite{frittelli2019preserving}. 
Following \cite{dautilia2020matrix}, the fully discrete scheme \eqref{semilinear_heat_equation_fully_discrete} can be further accelerated through an a-priori LU factorisation of the coefficient matrices (in combination with \texttt{symamd} reordering of such matrix to further increase sparsity). In matrix-oriented form, \eqref{semilinear_heat_equation_fully_discrete_pk} and \eqref{semilinear_heat_equation_fully_discrete} become
\begin{align}
\label{semilinear_heat_equation_fully_discrete_sylvester_pk}
\widetilde{\mathcal{L}}(U^{(k+1)}) = \widetilde{\mathcal{R}}(U^{(k)}), \qquad k=0,\dots, N_T,\\
\label{semilinear_heat_equation_fully_discrete_sylvester}
\widehat{\mathcal{L}}(U^{(k+1)}) = \widehat{\mathcal{R}}(U^{(k)}), \qquad k=0,\dots, N_T,
\end{align}
respectively, where $U^{(k)}$ is such that $\texttt{vec}(U^{(k)}) = \boldxi^{(k)}$ and
\begin{align}
\label{semilinear_heat_equation_fully_discrete_L_pk}
&\widetilde{\mathcal{L}}(U) := M_3UM +d_u\tau \Big(M_1UA + B_1UM + M_2UB_2 - C_2UC_1 - C_2^TUC_1^T\Big);\\
\label{semilinear_heat_equation_fully_discrete_L}
\begin{split}
&\widehat{\mathcal{L}}(U) := M_0D_1^2UM_0 +d_u\tau \Big((M_0UA + D_1A_1UM_0 + M_0D_2^2UA_2)\\
&\hspace*{10mm} -D_1(CD_2UCD_3 + (CD_2)^TU(CD_3)^T)\Big);
\end{split}\\
\label{semilinear_heat_equation_fully_discrete_R_pk}
&\widetilde{\mathcal{R}}(U) := M_3(U + \tau f(U))M;\\
\label{semilinear_heat_equation_fully_discrete_R}
&\widehat{\mathcal{R}}(U) := M_0D_1^2(U + \tau f(U))M_0.
\end{align}
Each iteration of the fully discrete scheme \eqref{semilinear_heat_equation_fully_discrete_sylvester_pk} or \eqref{semilinear_heat_equation_fully_discrete_sylvester} is a multiterm Sylvester equation that can be solved through the matrix-oriented PCG method \eqref{laplace_xnormal_pcg}-\eqref{laplace_xnormal_pcg_iteration} where $\mathcal{L}$ and $\mathcal{R}$ are chosen accordingly. For the preconditioning operator $\mathcal{P}$, observe in \eqref{semilinear_heat_equation_fully_discrete_sylvester_pk} that, in the limit $\frac{\tau}{N^2} \rightarrow 0$, the operator $\mathcal{L}(U)$ reduces to $M_3UM$. Hence, $\mathcal{P}(U) := M_3UM$ would make for a reasonable preconditioner. However, we experimentally found that by adding suitable corrections we obtain a more accurate preconditioner that retains a single-term form and it is given by:
\begin{align}
\label{semilinear_heat_equation_fully_discrete_P}
&\widetilde{\mathcal{P}}^{-1}(U) = (M_3+d_u\tau B_1)^{-1}U(M+d_u\tau(A+B_2))^{-1};\\
\label{semilinear_heat_equation_fully_discrete_P_lumped}
&\widehat{\mathcal{P}}^{-1}(U) = (M_0D_1+d_u\tau A_1)^{-1}U(M_0+d_u\tau(A+A_2))^{-1};
\end{align}
We use the following stop criterion: at each timestep, we stop the iterations of \eqref{laplace_xnormal_pcg_iteration} when the truncated solution $U^{(s)}$ fulfils $\|R^{(s)}\| \leq \tau\|R^{(0)}\|$, where $R^{(s)}$ is the residual defined by $R^{(s)} : = \mathcal{L}(U^{(s)}) - \mathcal{R}(U^{(s)})$, with the corresponding operators as defined in \eqref{semilinear_heat_equation_fully_discrete_L_pk}-\eqref{semilinear_heat_equation_fully_discrete_R}.
\footnote{This stopping criterion guarantees optimal convergence in space and time, and is justified as follows. At each timestep, the initial guess $U^{(0)}$ is an $O(\tau)$-accurate approximation of the exact solution $\widehat{U}$ of \eqref{semilinear_heat_equation_fully_discrete_sylvester} which in turn contains an $O(\tau^2)$ discretisation error, since IMEX-Euler is first-order accurate. Hence, to preserve the accuracy of the method, the truncated solution $U^{(s)}$ of \eqref{laplace_xnormal_pcg}-\eqref{laplace_xnormal_pcg_iteration} must be an $O(\tau^2)$-accurate approximation of $\widehat{U}$ as well. Consequently, $U^{(s)}$ must approximate $\widehat{U}$ better than $U^{(0)}$ by $O(\tau)$ times. In terms of residuals, $\|R^{(s)}\|$ must be $O(\tau)$ times $\|R^{(0)}\|$ in \eqref{laplace_xnormal_pcg}-\eqref{laplace_xnormal_pcg_iteration}.}

\subsection{Numerical Example: Semilinear heat equation on $x$-normal domain}
We consider the following heat equation with zero Dirichlet boundary conditions on the the cap-shaped domain $\Omega^S$ defined in \eqref{cap-shaped-domain} :
\begin{equation}
\label{example_heat_cap-shaped-domain_dirichlet}
\begin{cases}
&u_t -d_u\Delta u = f(x,y,t), \qquad (x,y) \in \Omega^S, \quad t\in [0,1];\\
&u(x,y,t) = 0, \qquad (x,y) \in  \partial \Omega^S, \quad t\in [0,1];\\
&u(x,y,0) = y(y - 1)\left(- \frac{y^2}{2} + x + 1\right)\left(\frac{y^2}{2} + x - 1\right), \qquad (x,y) \in \Omega^S,
\end{cases}
\end{equation}
were $d_u = 0.1$ and $f(x,y,t)$ is chosen in such a way that the exact solution is $u(x,y,t) = u(x,y,0)\exp(t)$, we omit the cumbersome expression of such $f(x,y,t)$. We consider both $\mathbb{P}_k$, $k=1,2,3,4$, and lumped $\mathbb{P}_1$ elements. Also in this case, we solve the vector formulations \eqref{semilinear_heat_equation_fully_discrete_pk} (for $\mathbb{P}_k$ elements) and \eqref{semilinear_heat_equation_fully_discrete} (for lumped $\mathbb{P}_1$ elements) via the MATLAB direct solver \texttt{mldivide}, with only one preliminary LU-decomposition. We compare these results with  the matrix-oriented PCG approach \eqref{semilinear_heat_equation_fully_discrete_sylvester} with preconditioner \eqref{semilinear_heat_equation_fully_discrete_P}.  We present a Test 1, to study the convergence of the two approaches and a Test 2 to highlight the computational advantages in time of the MO-PCG approach.

\begin{table}[t!]
\caption{\textbf{Test 1} - Semilinear heat equation \eqref{example_heat_cap-shaped-domain_dirichlet} on the cap-shaped domain:   $\mathbb{P}_k$ finite elements, $k=1,2,3,4$ and lumped $\mathbb{P}_1$ elements. We show the convergence rates in space and time obtained by solving the vector formulations \eqref{semilinear_heat_equation_fully_discrete_pk} {and \eqref{semilinear_heat_equation_fully_discrete}} by the direct method and the corresponding MO-FEM by the MO-PCG method \eqref{semilinear_heat_equation_fully_discrete_sylvester_pk}. The convergence rates are optimal for all $k=1,2,3,4$.}
\begin{center}
\begin{tabular}{ l c c c }
$k$ & \thead{Convergence rate \\ Vector method} & \thead{Convergence rate\\ Matrix PCG}\\ 
 \hline
 $1$ lumped & 1.9953 & 2.0509\\
 $1$  & 1.9945 & 1.9946\\  
 $2$  & 2.9945 & 2.9947\\
 $3$  & 3.9941 & 3.9943\\
 $4$  & 4.9939 & 4.9942\\
\end{tabular}
\end{center}
\label{tab:parabolic_convergence_pk}
\end{table}

In \textbf{Test 1}, for both methods we consider the second and third mesh $\Omega_i^S$, $i=1,\dots,2$ of Experiment \ref{sec:experiment_elliptic_cap-shaped-domain_curved}  for $N=48, 96$. Correspondingly, for each $i=1,2$ and $k=1,2,3,4$ we choose $\tau_{i,k} = 0.01 \cdot 2^{(k+1)(i-1)}$ and $\texttt{tol} = 1$. This choice of timesteps allows to highlight optimal convergence in $L^2(\Omega^S)$ norm (i.e. $(k+1)$-th order in space and first order in time). We have confined this test to these $N$ values because the timestep $\tau_{i,k}$ would become too small for larger values of $N$ when $k>1$. The results are shown in Table \ref{tab:parabolic_convergence_pk}.

In \textbf{Test 2}, we consider $N =480, 960, 1920$, fixed $\tau = 1$e-2. This test is more representative of typical user-case scenarios where high spatial resolution is required. The obtained results indicate a significant advantage of MO-PCG and are shown in Table \ref{tab:parabolic_times_pk} for all $k$ and for all $N$. At each timestep,  MO-PCG converges with just one iteration, except with lumped $\mathbb{P}_1$ elements, where up to two PCG iterations per timestep are required.

\begin{table}
\caption{\textbf{Test 2} - Semilinear heat equation \eqref{example_heat_cap-shaped-domain_dirichlet} on the cap-shaped domain \eqref{cap-shaped-domain}, solved through $\mathbb{P}_k$ finite elements, $k=1,2,3,4$ and lumped $\mathbb{P}_1$ elements. For each meshsize $N$, we apply the IMEX Euler method with timestep $\tau=0.01$. In all cases,  the MO-PCG approach is quicker than the vector (direct) approach, especially for $k=1$ and $k=4$. The gap increases with $N$, as shown by comparing the time ratios. MO-PCG always converges with one single iteration, except with lumped $\mathbb{P}_1$ elements, where two iterations are required.}
\begin{center}
\begin{tabular}{l l c c c c}
$N$ & $k$ & \thead{Time (s) \\ Vector method} & \thead{Time (s)\\ Matrix PCG} & \thead{Time ratio \\ (Vector/PCG)} & \thead{Iterations\\ PCG}\\ 
 \hline
\multirow{ 5}{*}{$480$} & $1$ lumped & 16.79 & 4.089 & 4.106 & 2\\
 & $1$  &  10.10 & 4.313 & 4.228 & 1\\ 
 & $2$  &  10.21 & 4.571 & 2.342 & 1\\
 & $3$  &  13.42 & 4.886 & 2.233 & 1\\
 & $4$  &  17.30 & 4.094 & 2.746 & 1\\ 
 \hline
\multirow{ 5}{*}{$960$} & $1$ lumped & 129.1 & 18.58 & 6.949& 2\\
 & $1$  &  131.7 & 16.41 & 8.023 & 1\\ 
 & $2$  &  54.98 & 19.10 & 2.879 & 1\\
 & $3$  &  58.79 & 20.25 & 2.902 & 1\\
 & $4$  &  89.18 & 21.62 & 4.126 & 1\\ 
 \hline
\multirow{ 5}{*}{$1920$} & $1$ lumped & 1269 & 77.50 & 16.37& 2\\
 & $1$  &  4337 & 74.35 & 58.33 & 1\\ 
 & $2$  &  497.9 & 87.01 & 5.722 & 1\\
 & $3$  &  473.0 & 81.17 & 5.828 & 1\\
 & $4$  &  3311 & 87.27 & 37.93 & 1\\ 
\end{tabular}
\end{center}
\label{tab:parabolic_times_pk}
\end{table}

\section{Applications to pattern formation in battery modeling}
\label{sec:turing_patterns}
We now consider the following reaction-diffusion model in two variables $\eta:\Omega \times [0,T]\rightarrow \mathbb{R}$ and $\theta:\Omega \times [0,T] \rightarrow [0,1]$, endowed with zero Neumann boundary conditions, on an arbitrary compact domain $\Omega \subset \mathbb{R}^2$:
\begin{equation}
\label{dib_model}
\begin{cases}
\eta_t - \Delta \eta = \rho f(\eta,\theta), \qquad (x,y,t) \in \Omega\times [0,T];\\
\theta_t - d_\theta\Delta \theta = \rho g(\eta,\theta), \qquad (x,y,t) \in \Omega\times [0,T];\\
\nabla \eta \cdot \boldn = \nabla \theta \cdot \boldn = 0, \qquad (x,y,t) \in \partial \Omega \times [0,T];\\
\eta(x,y,0) = \eta_0(x,y), \quad \theta(x,y,0) = \theta_0(x,y), \qquad (x,y) \in \Omega,
\end{cases}
\end{equation}
where $d_\theta>0$ is the diffusion coefficient, $\rho > 0$ is a space-time rescaling factor and the kinetics are
\begin{align}
\label{dib_kinetics_f}
&f(\eta,\theta) := A_1 (1-\theta)\eta - A_2\eta^3 - B(\theta-\alpha);\\
\label{dib_kinetics_g}
&g(\eta, \theta) := C(1+k_2\eta)(1-\theta)[1-\gamma(1-\theta)] - D\theta(1+\gamma\theta)(1+k_3\eta),
\end{align}
with $\alpha, \gamma, A_1, A_2, B, C, D, k_2, k_3$ positive parameters. 

The PDE system \eqref{dib_model}-\eqref{dib_kinetics_g} is known as DIB model and has been introduced for the first time in \cite{bozzini2013spatio} to describe electrodeposition processes. Under suitable choices of the parameters and of the domain $\Omega$, this model was shown to possess a variety of spatially structured solutions, known as Turing patterns, see for example \cite{lacitignola2014spatio, lacitignola2015spatio}. An interesting application in battery modeling is reported in \cite{lacitignola2017turing, lacitignola2019spiral}. Turing patterns are obtained as stationary solutions of \eqref{dib_model}-\eqref{dib_kinetics_g} and then their numerical approximation requires highly spatial accuracy for longtime integration, this motivates the development of efficient solvers. In this direction, a first work based on matrix oriented formulation of \eqref{dib_model}-\eqref{dib_kinetics_g} is \cite{dautilia2020matrix} where finite differences and several time solvers have been proposed on square domains. In \cite{dautilia2020matrix}, the Sylvester matrix equations obtained at each time step have been approximated by the \emph{reduced approach}, similar to the one in Section 4. For example the IMEX Euler yielded the \emph{rEuler} method, that revealed much more efficient than its classical vector approach.
On the other hand, domain geometry was also proven to play an important role in pattern selection, as shown also in \cite{lacitignola2017turing, lacitignola2019spiral}, for this reason efficient solvers that can be applied on domains as general as possible are need.
Towards this aim, here we propose the matrix oriented FEM spatial approximation and the MO-PCG approach presented in Setion 6 with preconditioner \eqref{semilinear_heat_equation_fully_discrete_P} to deal in particular with some x-normal domains.
We will present two kind of simulations, first on the cap shaped domain introduced in \eqref{cap-shaped-domain} and then on the jar-domain shown in Fig. \ref{fig:jar_shaped_domain} that correspond to the curvilinear cylinder in \ref{fig:wrapper_jar}. In both cases we will consider domains of increasing \emph{effective domain size} given by ${\cal A}= \rho | \Omega |$, where $|\Omega|$ is the area of the domain in \eqref{dib_model}. In fact, as shown in \cite{lacitignola2017turing, lacitignola2019spiral} there exists a sufficiently large  $\cal{A}^*$ such that for  $\cal{A} >\cal{A}^*$  the \emph{intrinsic} Turing pattern corresponding to the fixed model parameters arises (see also \cite{lacitignola2017turing} for more details), otherwise only a portion of it can be approximated, giving rise to doubts about its classification. 

To solve on domains of large sizes, we exploit the meaning and the role of the parameter $\rho$ in \eqref{dib_model}, as follows. By introducing new variables $(\widetilde{x}, \widetilde{y}) = \sqrt{\rho}(x,y)$, the chain rule yields
\begin{align}
\label{dib_eta_transformation}
\frac{\partial \eta}{\partial t} = \rho\frac{\partial \eta}{\partial\widetilde{t}}, \qquad \nabla_{(x,y)} \eta = \sqrt{\rho}\nabla_{(\widetilde{x}, \widetilde{y})} \eta, \qquad \Delta_{(x,y)}\eta = \rho\Delta_{(\widetilde{x}, \widetilde{y})} \eta;\\
\label{dib_theta_transformation}
\frac{\partial \theta}{\partial t} = \rho\frac{\partial \theta}{\partial\widetilde{t}}, \qquad \nabla_{(x,y)} \theta = \sqrt{\rho}\nabla_{(\widetilde{x}, \widetilde{y})} \theta, \qquad \Delta_{(x,y)}\theta = \rho\Delta_{(\widetilde{x}, \widetilde{y})} \theta;
\end{align}
Hence, we define $\Omega_\rho := \sqrt{\rho}\Omega$ and $T_\rho := \rho T$.
Hence, $\rho$ acts as a rescaling parameter in space and time. 

In the following simulations we always solve 110 in the reference domanis $\Omega$ cap shaped e jar shaped with final time time $T$ and timespet $\tau$ such that $T_\rho = 300$ and $\tau_{\rho} = 5e-3$ which guarantees the stability of the IMEX Euler method.  The corresponding numerical solutions will be plotted in the rescaled domain $\Omega_\rho$ for $T_\rho = 300$.
In all the experiments fix the following model parameters:
\begin{equation}
\label{dib_parameters_mixed}
\alpha = 0.5, \gamma = 0.2, A_1 = 10,  D = 3.2727, k_2 = 2.5, k_3 = 1.5, d_\theta=20.
\end{equation}
The initial data are given by $\theta_0(x,y) = \theta_e + 10^{-4}rand(x,y)$ and $\eta_0(x,y) = \eta_e + 10^{-4}rand(x,y)$ and are small spatially random perturbations of the homogeneous equilibrium $(\eta_e, \theta_e) := (0, 0.5)$.

\subsection{{Cap-shaped domain}}
In this example, we consider the cap-shaped domain \eqref{cap-shaped-domain} and we choose $A_2 = 30, B = 25, C = 7$,  that, according to the segmentation results in \cite{sgura2019parameter}, can yield mixed spots-worms Turing patterns at the steady state. 

We solve the PDE RDS system with $\mathbb{P}_k$ elements ($k=1,2,3,4$) and lumped $\mathbb{P}_1$ elements in space and by the IMEX Euler, as in the previous section. We compare the vector approach solving the sequence of linear systems in \eqref{semilinear_heat_equation_fully_discrete_pk}-\eqref{semilinear_heat_equation_fully_discrete} by the direct method (that we will call "vector method") with the MO-PCG approach \eqref{semilinear_heat_equation_fully_discrete_sylvester}, solving the multiterm Sylvester equations arising at each time step for this choice of the domain.

\begin{figure}[t!]
\includegraphics[scale=0.3]{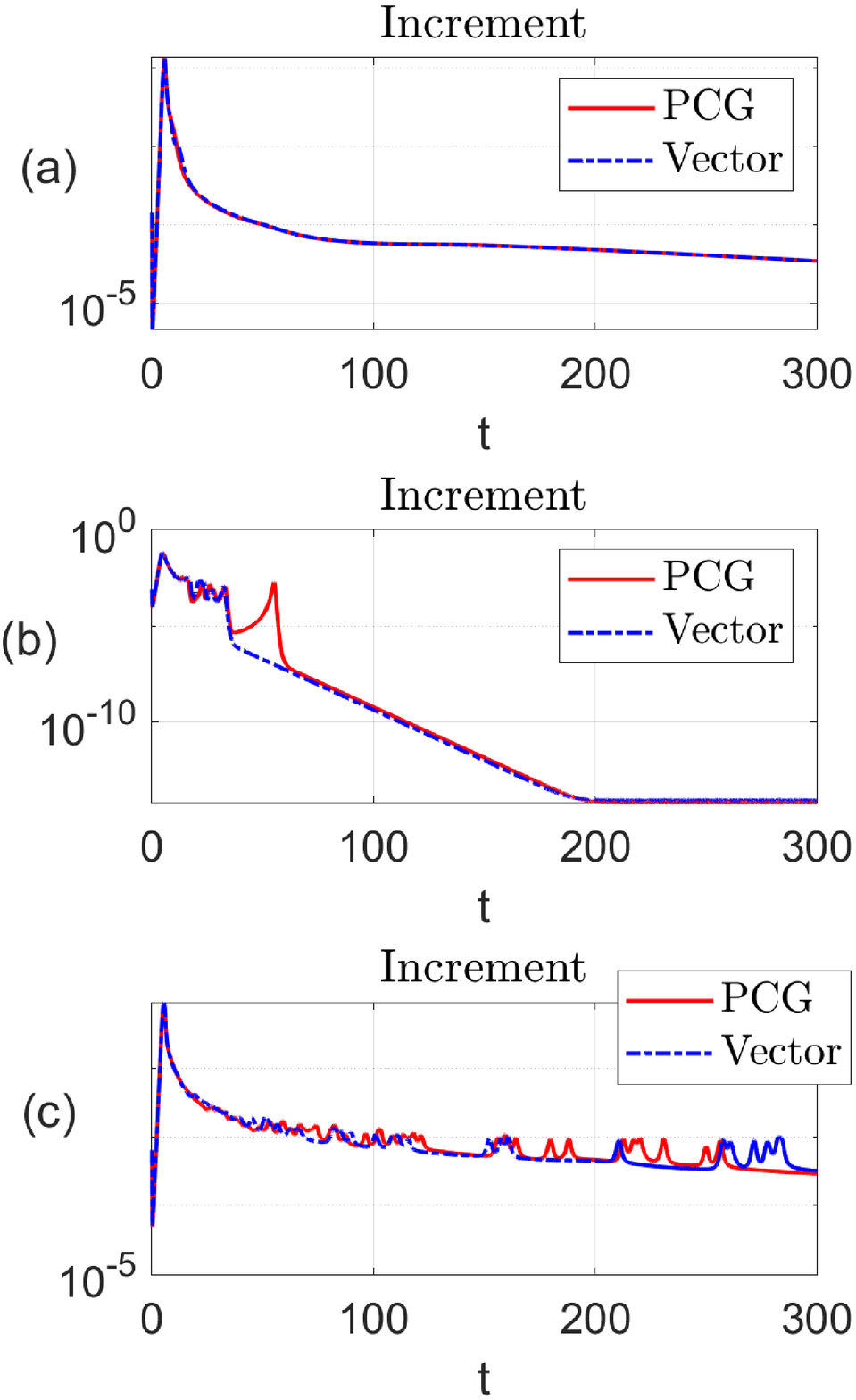}
\includegraphics[scale=0.3]{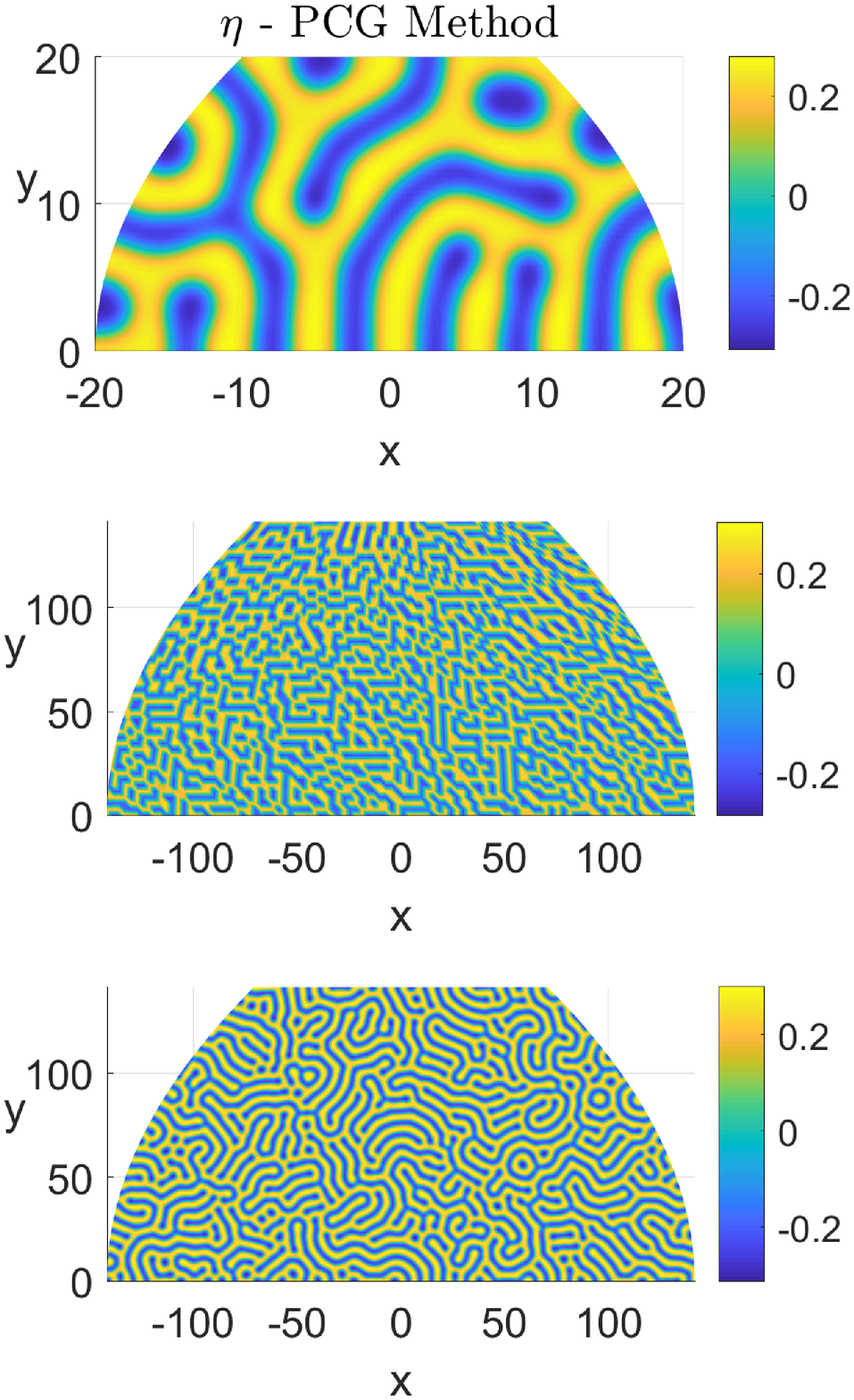}
\includegraphics[scale=0.3]{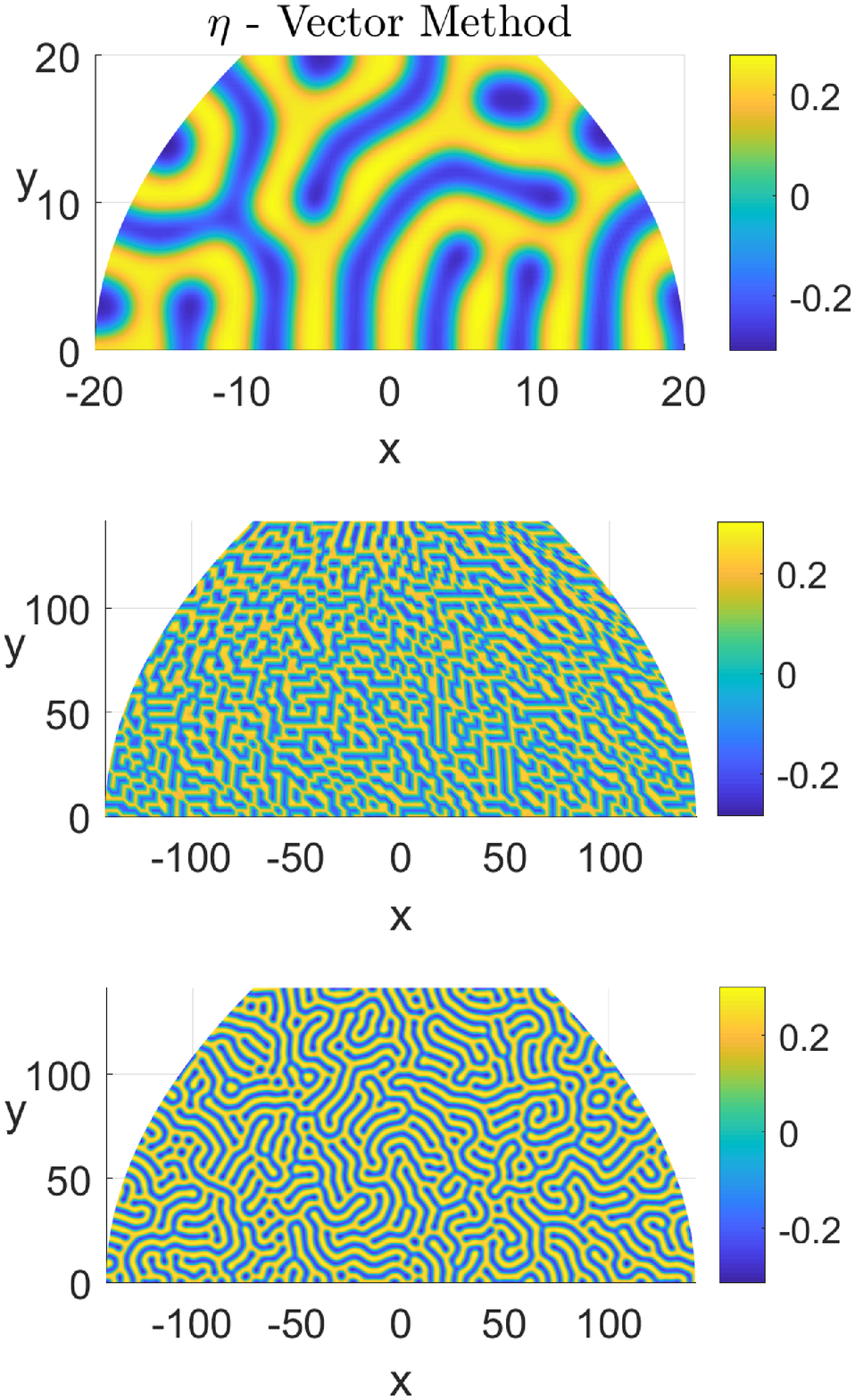}
\caption{DIB model \eqref{dib_model}-\eqref{dib_kinetics_g} on the cap-shaped domain, lumped $\mathbb{P}_1$ solutions.  The solutions are plotted on the rescaled domain $\Omega_\rho$. Each row corresponds to the $(\mathcal{A}, \rho,  N_x,  N_y)$ combination in Table \ref{tab:dib_times}. For $\rho = 400$ (a) only few structures arise in the pattern. For $\rho = 20000$: in (b) $N_x=200,N_y=100$ are not sufficient to resolve the pattern structure, that instead is well-resolved for $N_x = 400$ and $N_y = 200$ in (c). Smaller values of $N_x$ and $N_y$ do not provide sufficient spatial approximation on the larger domain and the pattern appear grainy (b). The patterns obtained by the MO-PCG method and the vector (direct) method are very similar, they are stationary solutions obtained at $T_\rho=300$, are shown by the increment dynamics (left subplots). To compare the execution times by the lumped $\mathbb{P}_1$ approximation see Table 4.}
\label{fig:spot_laby}
\end{figure}

We solve the DIB model \eqref{dib_model} with different combinations of $\rho$, $N_x$ and $N_y$ as listed in Table \ref{tab:dib_times}, that is for domains of larger area $\mathcal{A} = |\Omega_\rho|$. 
In all the computations, we discretise the $x$ dimension with $N_x$ nodes and the $y$ dimension with $N_y$ nodes, with $N_x = 2N_y$, which reflects the aspect ratio of the domain. 
In  Fig. \ref{fig:spot_laby}, for each simulation, we report the final patterns obtained by the vector approach and by the MO-PCG, together with the respective increments $\|\eta^{(k+1)} - \eta^{(k)}\|_F$ as a function of time, with $\|\cdot\|_F$ Frobenius norm. If such increment decreases over time and tends to an almost small stationary value, we deduce that the numerical solution is converging to a steady state. 

For increasing values of the \emph{effective domain size} $\cal A$ the solution morphology changes and a pattern with more structures is attained, as shows in Fig.6,(a)--(c) corresponding to the values (a)--(c) in Table 4, respectively.
In simulation (a) a good pattern is attained by both methods, but its morphology is not completely expressed. The vector and the matrix approach seem to be equivalent in this case also in terms of computational times (see Table 4).
To capture the true Turing morphology a larger domain $\Omega$ and a \emph{sufficiently} fine mesh is required, otherwise \emph{phantom patterns} could be obtained. This is exactly what happens for the simulation in case (b), corresponding to the second row of both Fig.6 and Table 4, where the same mesh of case (a) yields a ``pixelated'' pattern. Hence, for the same domain, in simulation (c) a finer mesh is used and both methods are finally able to attain a ``complete'' pattern. Note that, in vector form at each time step we solve a problem of dimension $N_x\cdot N_y =200 \cdot 400=8 \cdot 10^4$, by the MO-PCG instead we solve a sequence of $N_t=6\cdot 10^4$ rectangular multiterm Sylvester equations of size $ 400 \times 200$. 
Moreover, this example shows that the matrix-oriented PCG algorithm \eqref{laplace_xnormal_pcg}-\eqref{laplace_xnormal_pcg_iteration} can successfully solve rectangular Sylvester equations.
As shown in Table \ref{tab:dib_times}, the time ratios indicate that the matrix PCG approach \eqref{semilinear_heat_equation_fully_discrete_sylvester_pk}-\eqref{semilinear_heat_equation_fully_discrete_sylvester} tends to become quicker than the vector-direct approach \eqref{semilinear_heat_equation_fully_discrete} with significant advantage only for $\mathbb{P}_4$ and lumped and no lumped $\mathbb{P}_1$ elements. For $ k=2,3$, we guess that a different preconditioner could improve the results of the MO-PCG method.

\begin{table}
\caption{ DIB model \eqref{dib_model}-\eqref{dib_kinetics_g} with parameters \eqref{dib_parameters_mixed} on the cap-shaped domain \eqref{cap-shaped-domain}. Parameters for the simulations (a)-(b)-(c) reported in Figure 6 and performance comparison between the vector approach \eqref{semilinear_heat_equation_fully_discrete_pk}-\eqref{semilinear_heat_equation_fully_discrete} based on LU decomposition and the matrix PCG method \eqref{semilinear_heat_equation_fully_discrete_sylvester}. 
When the time ratio $r_t$ is close to 1, the methods take approximately the same time; $r_t > 1$ indicates that MO-PCG is quicker. The MO-PCG is less expensive as $N_x$ and $N_y$ increase. Simulation in (b) yields a \emph{phantom pattern}. The last column shows the amount of iterations required by MO-PCG for each of the two PDEs of the model.}
\begin{center}
\begin{tabular}{l | c c l l || l c c c c}
& $\cal{A}$ & $\rho$ & $N_x$ & $N_y$ & $k$ & \thead{Time (s) \\ Vector method} & \thead{Time (s)\\ Matrix PCG} & \thead{Time ratio\\ $r_t$} & \thead{Iterations\\ PCG $(\eta,\theta)$}\\ 
 \hline
\multirow{ 5}{*}{(a)  }  & \multirow{ 5}{*}{$2000/3$}  & \multirow{ 5}{*}{$400$} & \multirow{ 5}{*}{$100$} & \multirow{ 5}{*}{$50$} & $1$ lumped & 57.18 & 93.76 & 0.6099 & (1,6)\\
 & & & & & $1$   & 62.76 &  94.60& 0.6635 & (2,3)\\ 
  & & & & & $2$  & 75.60 &  182.0 & 0.4159 & (2,4)\\
  & &  & & & $3$  & 226.6 &  193.4  & 1.172 & (2,3)\\
  & &  & & & $4$  &  270.2 & 207.7 & 1.301 & (2,3)\\ 
 \hline
\multirow{ 5}{*}{(b)  }  & \multirow{ 5}{*}{$10000/3$}  & \multirow{ 5}{*}{$20000$} & \multirow{ 5}{*}{$100$} & \multirow{ 5}{*}{$50$} & $1$ lumped & 58.08 & 47.53 & 1.222 & (1,2)\\
  & & & & & $1$  & 69.72 &  70.89  & 0.9835 & (1,3)\\ 
  & & & & & $2$  & 75.50 &  118.1  & 0.6394 & (1,3)\\
  & & & & & $3$  & 82.79 &  142.6  & 0.5804 & (1,3)\\
  & & & & & $4$  & 273.7 &  156.5  & 1.749 & (1,3)\\
 \hline
\multirow{ 5}{*}{(c)  }  & \multirow{ 5}{*}{$10000/3$}  &  \multirow{ 5}{*}{$20000$} & \multirow{ 5}{*}{$400$} & \multirow{ 5}{*}{$200$} & $1$ lumped & 1590 & 879.8  & 1.807 & (1,4)\\
  & & & & & $1$  & 1623 &  1249  & 1.299 & (2,3)\\ 
  & & & & & $2$  & 1809 &  1934  & 0.9352 & (1,3)\\
  & & & & & $3$  & 1973  &  2143 & 0.9209 & (2,3)\\
  & & & & & $4$  & 6323 &  2266 & 2.791 & (2,3)\\ 
\end{tabular}
\end{center}
\label{tab:dib_times}
\end{table}

\subsection{{Jar-shaped domain}}
\label{sec:jar_shaped_domain}
\begin{figure}[t!]
\includegraphics[scale=0.3]{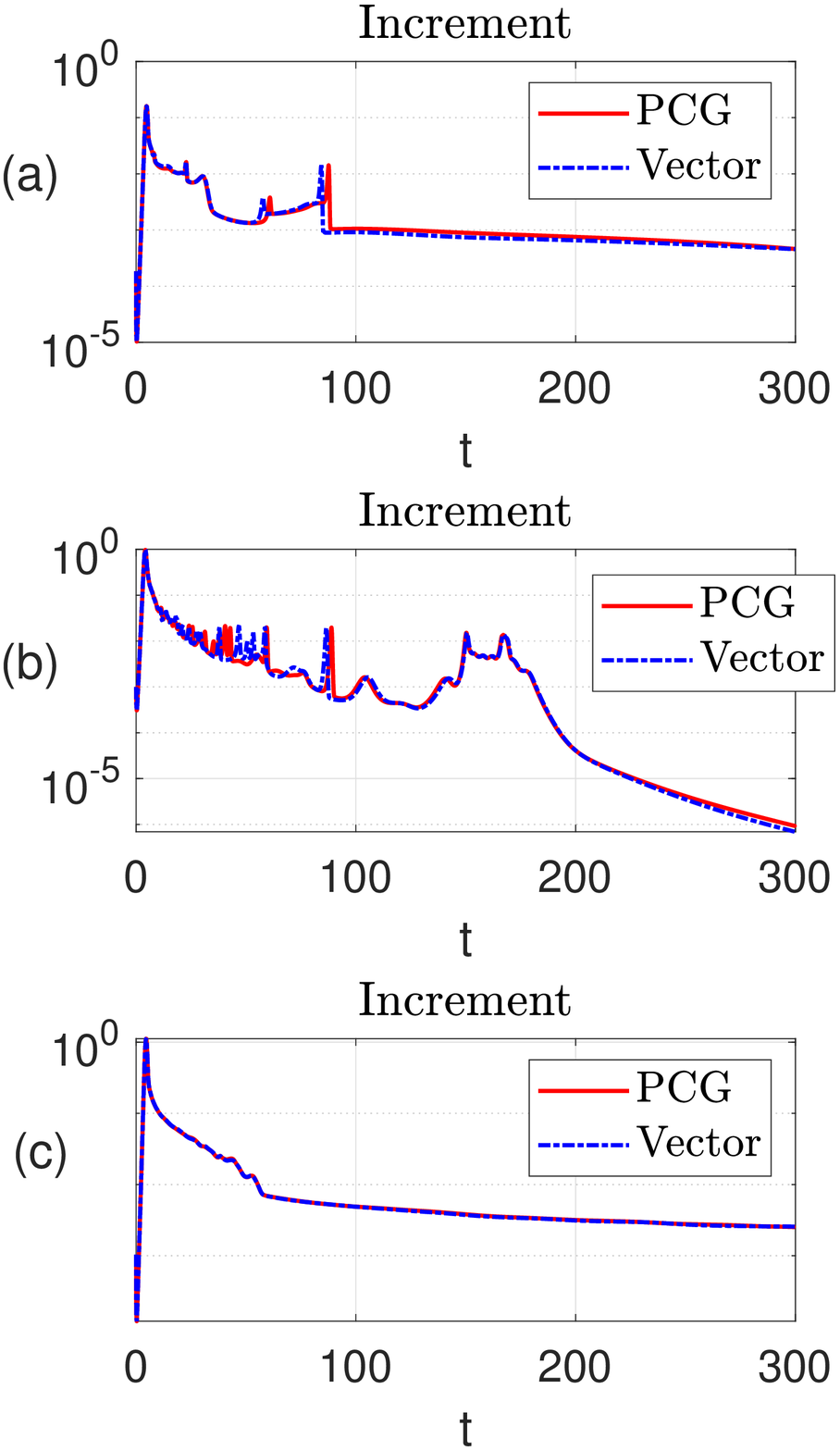}
\includegraphics[scale=0.3]{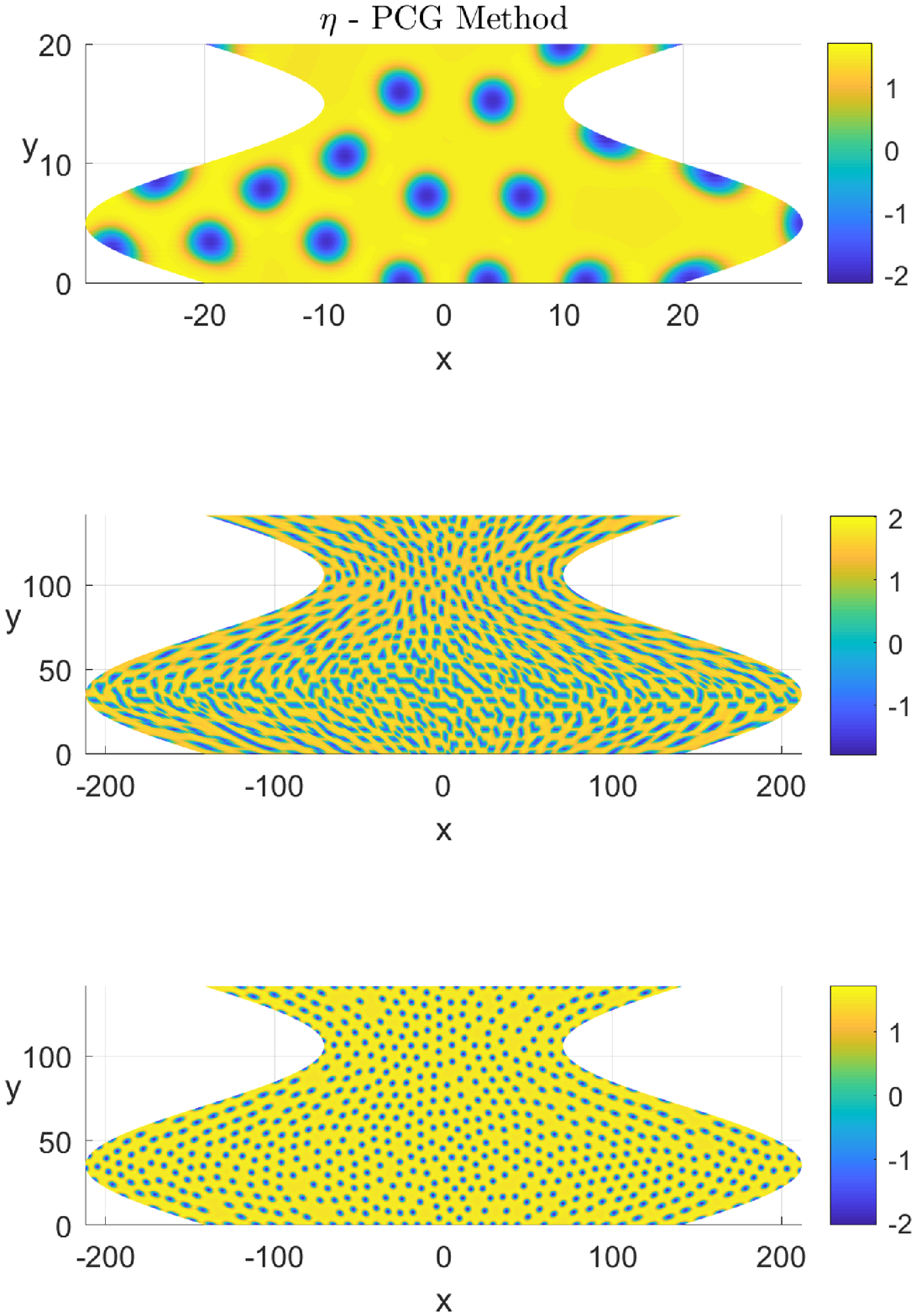}
\includegraphics[scale=0.3]{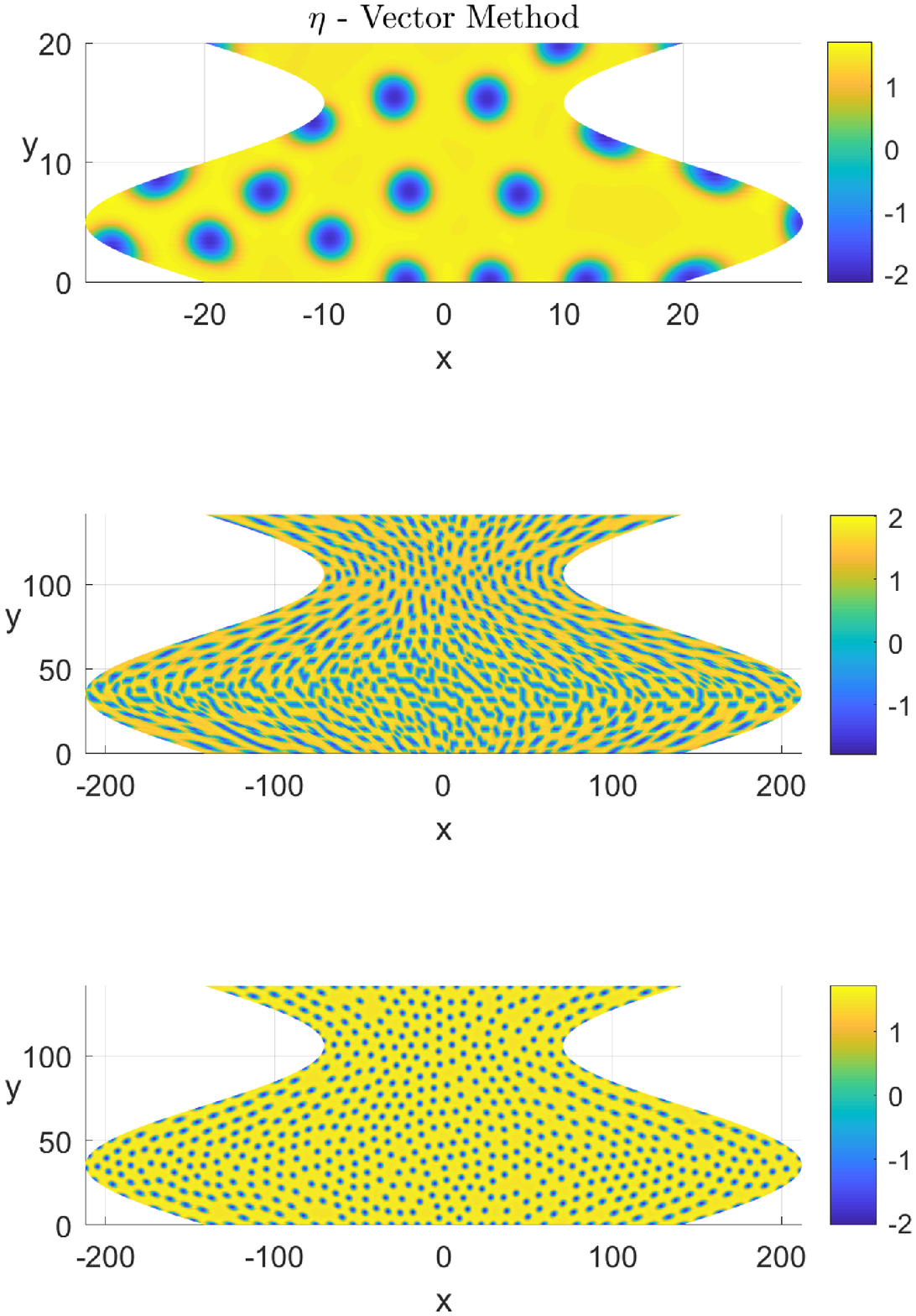}
\caption{DIB model \eqref{dib_model}-\eqref{dib_kinetics_g} on the jar-shaped domain, $\mathbb{P}_1$ solutions. Values of  $(\rho, N_x, N_y)$ are given in Table \ref{tab:dib_times_jar}.
In case (a) on the smallest domain only few holes arise in the pattern. In case (b), for $\rho = 20000$ on the larger domain the mesh is too coarse and a \emph{phantom pattern } arises. In case (c), the intrinsic Turing pattern is well-resolved for $N_x = 600$ and $N_y = 200$. MO-PCG and vector solutions are very similar, but the MO approach converges in significant less time (see Table 5).}
\label{fig:holes}
\end{figure}

Thanks to the results in the previous test, here we solve the model only with $\mathbb{P}_1$ elements, both with and without lumping. To further explore the robustness of the matrix PCG approach w.r.t. \emph{domain complexity and mesh distortion}, we consider the \emph{jar-shaped domain} in Fig. \ref{fig:jar_shaped_domain}. We fix the parameters $ A_2 = 1, B = 30, C = 3$ for the DIB model which are known to produce Turing patterns with holes (also called reversed spots) \cite{sgura2019parameter} and again we solve for increasing $\mathcal{A}$ with different combinations of $\rho$, $N_x$ and $N_y$ as listed in Table \ref{tab:dib_times_jar}.  In all the computations, we consider $N_x = 3N_y$, which reflects the aspect ratio of the domain. The timestep and the final time are as in the previous test on the cap-shaped domain. We show the $\mathbb{P}_1$ solutions in Fig. \ref{fig:holes}.\\
As we can see in the figure,  in case (a) on the smallest domain only few holes arise in the pattern. In case (b) on the larger domain, more structures arise in the pattern, but the mesh is too coarse and a \emph{phantom pattern} arises. In case (c),  the intrinsic Turing pattern is well-resolved for $N_x = 600$ and $N_y = 200$.  The MO-PCG and the vector solutions are very similar, but the MO approach converges in significant less time (see Table 5).

We conclude by remarking that, since the jar-shaped domain in Fig. \ref{fig:jar_shaped_domain} can be transformed to the cylinder $\Gamma$ shown in Fig. \ref{fig:wrapper_jar}, then the solutions shown in Fig. \ref{fig:holes} can be interpreted, after the coordinate transformation \eqref{transformation_cylinder}, as solutions to the \emph{surface DIB model}, that is \eqref{dib_model}-\eqref{dib_kinetics_g} where the Laplace operator $\Delta$ is replaced by the Laplace-Beltrami operator $\Delta_\Gamma$ on the cylinder $\Gamma$. As an example, we report in Fig. 8 the solution in (a) wrapped on a curvilinear cylinder. The application of the model on a cylindrical surface can be of applicative interest as shown in \cite{bozzini2020morphological}, in which the authors consider the use of cylindrical Zn sponges as a means of limiting the shape change and dendrite formation issues in Zn-based rechargeable batteries.

\begin{table}[t!]
\caption{Reaction-diffusion model \eqref{dib_model}-\eqref{dib_kinetics_g} with parameters as in Section \ref{sec:jar_shaped_domain} on the jar-shaped domain in Fig. \ref{fig:jar_shaped_domain}: performance comparison between the vector method \eqref{semilinear_heat_equation_fully_discrete_pk}-\eqref{semilinear_heat_equation_fully_discrete} with LU decomposition and the matrix PCG method \eqref{semilinear_heat_equation_fully_discrete_sylvester}. The vector PCG approach tends to become quicker than the vector approach as $N_x$ and $N_y$ increase.}
\begin{center}
\begin{tabular}{l | l l l l || l c c c c}
\ & $\mathcal{A}$ & $\rho$ & $N_x$ & $N_y$ & $k$ & \thead{Time (s) \\ Vector method} & \thead{Time (s)\\ Matrix PCG} & Time ratio & \thead{Iterations\\ PCG $(\eta,\theta)$}\\ 
 \hline
\multirow{ 2}{*}{(a)} & \multirow{ 2}{*}{$800$} & \multirow{ 2}{*}{$400$} & \multirow{ 2}{*}{$150$} & \multirow{ 2}{*}{$50$} & $1$ lumped & 309.6 & 306.6 & 1.010 & (2,16)\\
& & & & & $1$  & 328.7 &  278.4 & 1.181 & (3,7)\\ 
 \hline
\multirow{ 2}{*}{(b)} & \multirow{ 2}{*}{$40000$} &\multirow{ 2}{*}{$20000$} & \multirow{ 2}{*}{$150$} & \multirow{ 2}{*}{$50$} & $1$ lumped & 98.54 & 67.17 & 1.4670 & (1,2)\\
& & & & & $1$ & 111.4  & 117.0  & 0.9528 & (1,2)\\ 
 \hline
\multirow{ 2}{*}{(c)} & \multirow{ 2}{*}{$40000$} &\multirow{ 2}{*}{$20000$} & \multirow{ 2}{*}{$600$} & \multirow{ 2}{*}{$200$} & $1$ lumped & 7995 & 2316.4  & 3.451 & (1,9)\\
& & & & & $1$ & 8067 & 3028  & 2.664 & (2,7)\\ 
\end{tabular}
\end{center}
\label{tab:dib_times_jar}
\end{table}

\begin{figure}
\begin{center}
\includegraphics[scale=0.5]{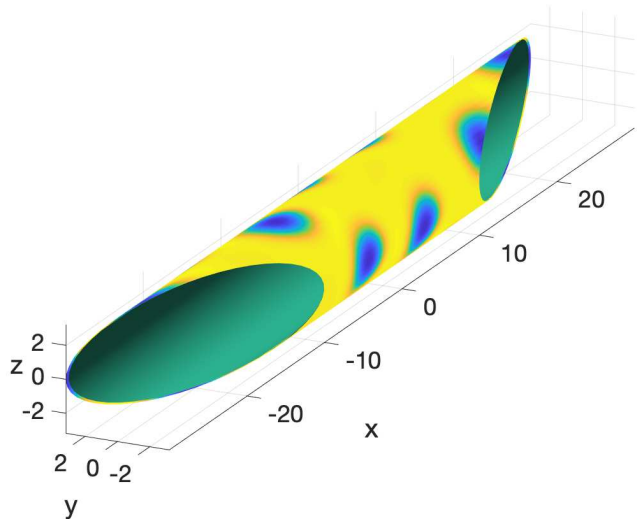}
\end{center}
\caption{Numerical solution of Fig. \ref{fig:holes} (a)  ($\rho=20000$, $N_x = 600$, $N_y = 200$) mapped onto the cylinder $\Gamma$ in Fig. \ref{fig:wrapper_jar}. It can be interpreted as the solution of the DIB \emph{surface reaction-diffusion model} on $\Gamma$.}
\end{figure}

\section{Conclusions}
\label{sec:conclusions}
{In this work we have provided a matrix-oriented formulation for Lagrangian finite elements of arbitrarily high order $k\in\mathbb{N}$ on $x$-normal domains. The proposed approach applies to both elliptic and parabolic PDE problems. The discrete problems take the form of a matrix equation (or a sequence of matrix Sylvester equations in the time-dependent case) of much smaller dimension that is mathematically equivalent to the much larger standard linear systems in Kronecker form. The proposed approach adopts a curvilinear structured mesh that eliminates geometric boundary error. Moreover, through a coordinate transformation, our approach applies also to special surface domains, namely cylinders with arbitrary curved boundaries.\\
On square domains, the discrete problems take the form of a generalised two-term Sylvester equation that we solve efficiently through a spectral approach for all $k$. In this sense, our work extends the findings in \cite{dautilia2020matrix}, based on classical finite differences, to the case of high order FEM in space. \\
On general $x$-normal domains, the discrete problem takes the form of (a sequence of) \emph{multiterm} Sylvester equations which we solve through a matrix-oriented PCG method with matrix-oriented preconditioner that is quick to evaluate thanks to is single-term form. On one hand, such solver is always quicker than the classical PCG in vector form. On the other hand, we show by several numerical tests that it is quicker than the optimised direct solver \texttt{mldivide} of MATLAB in the case of (i) time-dependent PDEs on general $x$-normal domains and (ii) elliptic PDEs on square domains. In terms of memory occupation, the matrix-oriented PCG method always improves on any direct or iterative solver that relies on the full storage of the Kronecker matrix, and the gap increases both with the number of gridpoints $N$ and the polynomial order $k$ of the method. The only case in which we could not find any speedup in the MO-PCG approach is that of elliptic problems on non-square domains with $\mathbb{P}_k$ elements, $k \neq 1$. This opens the quest for efficient preconditioners, which will be addressed in future studies.\\
Special consideration deserves the application to reaction-diffusion systems, where the simulation of fine-grained patterns requires high spatial resolution, which translates into computationally intense simulations both in time and memory. Our experiments for the approximation of Turing patterns arising in batteries as solution of the DIB morphochemical model provide encouraging results in this direction and justify the matrix approach in terms of execution times and storage. 
The best performance gains were observed with $\mathbb{P}_4$ and lumped $\mathbb{P}_1$ elements.  Also in this case, we believe that further performance gains can be found through the development of more efficient preconditioners and more efficient solvers for multiterm Sylvester equations, such as a truncated PCG \cite{shank2016efficient}. These aspects will be addressed in future studies.

\section*{Acknowledgments}
The work of MF was funded by Regione Puglia (Italy) through the research programme REFIN-Research for Innovation (protocol code 901D2CAA, project number UNISAL026).\\
The work of IS has been funded by the MIUR (Italian Ministry of Education, University and Research) project PRIN 2017, ``Mathematics of active materials: From mechanobiology to smart devices'', project no. 2017KL4EF3.\\
The work of MF and IS was performed under the auspices of GNCS-INdAM (Italian National Group of Scientific Computing).
\FloatBarrier
 
\bibliographystyle{plainurl}
\bibliography{bibliography}

\begin{thebibliography}{10}

\bibitem{antolin2015efficient}
P~Antolin, A~Buffa, F~Calabr\'{o}, M~Martinelli, and G~Sangalli.
\newblock Efficient matrix computation for tensor-product isogeometric
  analysis: The use of sum factorization.
\newblock {\em Computer Methods in Applied Mechanics and Engineering},
  285:817--828, 2015.
\newblock \href {https://doi.org/10.1016/j.cma.2014.12.013}
  {\path{doi:10.1016/j.cma.2014.12.013}}.

\bibitem{barreira2011surface}
R~Barreira, C~M Elliott, and A~Madzvamuse.
\newblock The surface finite element method for pattern formation on evolving
  biological surfaces.
\newblock {\em Journal of Mathematical Biology}, 63(6):1095--1119, 2011.
\newblock \href {https://doi.org/10.1007/s00285-011-0401-0}
  {\path{doi:10.1007/s00285-011-0401-0}}.

\bibitem{bartels1972solution}
R~H Bartels and G~W Stewart.
\newblock Solution of the matrix equation {AX + XB = C}.
\newblock {\em Communications of the ACM}, 15(9):820--826, 1972.
\newblock \href {https://doi.org/10.1145/361573.361582}
  {\path{doi:10.1145/361573.361582}}.

\bibitem{becherer2005classical}
D~Becherer, M~Schweizer, et~al.
\newblock Classical solutions to reaction--diffusion systems for hedging
  problems with interacting {It{\^o}} and point processes.
\newblock {\em Annals of Applied Probability}, 15(2):1111--1144, 2005.
\newblock \href {https://doi.org/10.1214/105051604000000846}
  {\path{doi:10.1214/105051604000000846}}.

\bibitem{bozzini2013spatio}
B~Bozzini, D~Lacitignola, and I~Sgura.
\newblock Spatio-temporal organization in alloy electrodeposition: a
  morphochemical mathematical model and its experimental validation.
\newblock {\em Journal of Solid State Electrochemistry}, 17(2):467--479, 2013.
\newblock \href {https://doi.org/10.1007/s10008-012-1945-7}
  {\path{doi:10.1007/s10008-012-1945-7}}.

\bibitem{bozzini2020morphological}
B~Bozzini, C~Mele, A~Veneziano, N~Sodini, G~Lanzafame, A~Taurino, and
  L~Mancini.
\newblock Morphological evolution of {Zn}-sponge electrodes monitored by in
  situ {X}-ray computed microtomography.
\newblock {\em ACS Applied Energy Materials}, 3(5):4931--4940, 2020.
\newblock \href {https://doi.org/10.1021/acsaem.0c00489}
  {\path{doi:10.1021/acsaem.0c00489}}.

\bibitem{chaplain2001spatio}
M~A~J Chaplain, M~Ganesh, and I~G Graham.
\newblock Spatio-temporal pattern formation on spherical surfaces: numerical
  simulation and application to solid tumour growth.
\newblock {\em Journal of Mathematical Biology}, 42(5):387--423, 2001.
\newblock \href {https://doi.org/10.1007/s002850000067}
  {\path{doi:10.1007/s002850000067}}.

\bibitem{chung2000discrete}
F~Chung and S-T Yau.
\newblock Discrete {Green}'s functions.
\newblock {\em Journal of Combinatorial Theory, Series A}, 91(1-2):191--214,
  2000.
\newblock \href {https://doi.org/10.1006/jcta.2000.3094}
  {\path{doi:10.1006/jcta.2000.3094}}.

\bibitem{dautilia2020matrix}
M~C D'Autilia, I~Sgura, and V~Simoncini.
\newblock Matrix-oriented discretization methods for reaction--diffusion
  {PDE}s: {Comparisons} and applications.
\newblock {\em Computers \& Mathematics with Applications}, 79(7):2067--2085,
  2020.
\newblock \href {https://doi.org/10.1016/j.camwa.2019.10.020}
  {\path{doi:10.1016/j.camwa.2019.10.020}}.

\bibitem{eilks2008numerical}
C~Eilks and C~M Elliott.
\newblock Numerical simulation of dealloying by surface dissolution via the
  evolving surface finite element method.
\newblock {\em Journal of Computational Physics}, 227(23):9727--9741, 2008.
\newblock \href {https://doi.org/10.1016/j.jcp.2008.07.023}
  {\path{doi:10.1016/j.jcp.2008.07.023}}.

\bibitem{elliott2010modeling}
C~M Elliott and B~Stinner.
\newblock Modeling and computation of two phase geometric biomembranes using
  surface finite elements.
\newblock {\em Journal of Computational Physics}, 229(18):6585--6612, 2010.
\newblock \href {https://doi.org/10.1016/j.jcp.2010.05.014}
  {\path{doi:10.1016/j.jcp.2010.05.014}}.

\bibitem{frittelli2019preserving}
M~Frittelli, A~Madzvamuse, I~Sgura, and C~Venkataraman.
\newblock Preserving invariance properties of reaction--diffusion systems on
  stationary surfaces.
\newblock {\em IMA Journal of Numerical Analysis}, 39(1):235--270, 2019.
\newblock \href {https://doi.org/10.1093/imanum/drx058}
  {\path{doi:10.1093/imanum/drx058}}.

\bibitem{golub1979hessenberg}
G~Golub, S~Nash, and C~Van~Loan.
\newblock A {Hessenberg-Schur} method for the problem {AX + XB = C}.
\newblock {\em IEEE Transactions on Automatic Control}, 24(6):909--913, 1979.
\newblock \href {https://doi.org/10.1109/tac.1979.1102170}
  {\path{doi:10.1109/tac.1979.1102170}}.

\bibitem{golub2013matrix}
G~H Golub and C~F Van~Loan.
\newblock {\em Matrix computations}, volume~3.
\newblock JHU press, 2013.

\bibitem{hao2020matrix}
Yue Hao and Valeria Simoncini.
\newblock Matrix equation solving of {PDEs} in polygonal domains using
  conformal mappings.
\newblock {\em Journal of Numerical Mathematics}, 0(0), nov 2020.
\newblock \href {https://doi.org/10.1515/jnma-2020-0035}
  {\path{doi:10.1515/jnma-2020-0035}}.

\bibitem{hughes2012finite}
T~J~R Hughes.
\newblock {\em The finite element method: linear static and dynamic finite
  element analysis}.
\newblock Courier Corporation, 2012.

\bibitem{jordan1965calculus}
C~Jordan and K~Jord{\'a}n.
\newblock {\em Calculus of finite differences}, volume~33.
\newblock American Mathematical Soc., 1965.

\bibitem{kansa1990multiquadrics}
E~J Kansa.
\newblock Multiquadrics—a scattered data approximation scheme with
  applications to computational fluid-dynamics -- {II} solutions to parabolic,
  hyperbolic and elliptic partial differential equations.
\newblock {\em Computers \& Mathematics with Applications}, 19(8-9):147--161,
  1990.
\newblock \href {https://doi.org/10.1016/0898-1221(90)90271-k}
  {\path{doi:10.1016/0898-1221(90)90271-k}}.

\bibitem{lacitignola2017turing}
D~Lacitignola, B~Bozzini, M~Frittelli, and I~Sgura.
\newblock Turing pattern formation on the sphere for a morphochemical
  reaction-diffusion model for electrodeposition.
\newblock {\em Communications in Nonlinear Science and Numerical Simulation},
  48:484--508, 2017.
\newblock \href {https://doi.org/10.1016/j.cnsns.2017.01.008}
  {\path{doi:10.1016/j.cnsns.2017.01.008}}.

\bibitem{lacitignola2014spatio}
D~Lacitignola, B~Bozzini, and I~Sgura.
\newblock Spatio-temporal organization in a morphochemical electrodeposition
  model: analysis and numerical simulation of spiral waves.
\newblock {\em Acta Applicandae Mathematicae}, 132(1):377--389, 2014.
\newblock \href {https://doi.org/10.1007/s10440-014-9910-3}
  {\path{doi:10.1007/s10440-014-9910-3}}.

\bibitem{lacitignola2015spatio}
D~Lacitignola, B~Bozzini, and I~Sgura.
\newblock Spatio-temporal organization in a morphochemical electrodeposition
  model: {Hopf} and {Turing} instabilities and their interplay.
\newblock {\em European Journal of Applied Mathematics}, 26(2):143--173, 2015.
\newblock \href {https://doi.org/10.1017/s0956792514000370}
  {\path{doi:10.1017/s0956792514000370}}.

\bibitem{lacitignola2019spiral}
D~Lacitignola, I~Sgura, B~Bozzini, T~Dobrovolska, and I~Krastev.
\newblock Spiral waves on the sphere for an alloy electrodeposition model.
\newblock {\em Communications in Nonlinear Science and Numerical Simulation},
  79:104930, 2019.
\newblock \href {https://doi.org/10.1016/j.cnsns.2019.104930}
  {\path{doi:10.1016/j.cnsns.2019.104930}}.

\bibitem{mantzaflaris2017low}
A~Mantzaflaris, B~J{\"u}ttler, B~N Khoromskij, and U~Langer.
\newblock Low rank tensor methods in {Galerkin}-based isogeometric analysis.
\newblock {\em Computer Methods in Applied Mechanics and Engineering},
  316:1062--1085, 2017.
\newblock \href {https://doi.org/10.1016/j.cma.2016.11.013}
  {\path{doi:10.1016/j.cma.2016.11.013}}.

\bibitem{nie1985lumped}
Y-Y Nie and V~Thom{\'e}e.
\newblock A lumped mass finite-element method with quadrature for a non-linear
  parabolic problem.
\newblock {\em IMA Journal of Numerical Analysis}, 5(4):371--396, 1985.
\newblock \href {https://doi.org/10.1093/imanum/5.4.371}
  {\path{doi:10.1093/imanum/5.4.371}}.

\bibitem{palitta2016matrix}
D~Palitta and V~Simoncini.
\newblock Matrix-equation-based strategies for convection--diffusion equations.
\newblock {\em BIT Numerical Mathematics}, 56(2):751--776, 2016.
\newblock \href {https://doi.org/10.1007/s10543-015-0575-8}
  {\path{doi:10.1007/s10543-015-0575-8}}.

\bibitem{powell2017efficient}
C~E Powell, D~Silvester, and V~Simoncini.
\newblock An efficient reduced basis solver for stochastic {Galerkin} matrix
  equations.
\newblock {\em SIAM Journal on Scientific Computing}, 39(1):A141--A163, 2017.
\newblock \href {https://doi.org/10.1137/15m1032399}
  {\path{doi:10.1137/15m1032399}}.

\bibitem{saad2003iterative}
Y~Saad.
\newblock {\em Iterative methods for sparse linear systems}.
\newblock SIAM, 2003.

\bibitem{sangalli2016isogeometric}
G~Sangalli and M~Tani.
\newblock Isogeometric preconditioners based on fast solvers for the
  {Sylvester} equation.
\newblock {\em SIAM Journal on Scientific Computing}, 38(6):A3644--A3671, 2016.
\newblock \href {https://doi.org/10.1137/16m1062788}
  {\path{doi:10.1137/16m1062788}}.

\bibitem{sgura2019parameter}
I~Sgura, A~S Lawless, and B~Bozzini.
\newblock Parameter estimation for a morphochemical reaction--diffusion model
  of electrochemical pattern formation.
\newblock {\em Inverse Problems in Science and Engineering}, 27(5):618--647,
  2019.
\newblock \href {https://doi.org/10.1080/17415977.2018.1490278}
  {\path{doi:10.1080/17415977.2018.1490278}}.

\bibitem{shank2016efficient}
S~D Shank, V~Simoncini, and D~B Szyld.
\newblock Efficient low-rank solution of generalized {Lyapunov} equations.
\newblock {\em Numerische Mathematik}, 134(2):327--342, 2016.
\newblock \href {https://doi.org/10.1007/s00211-015-0777-7}
  {\path{doi:10.1007/s00211-015-0777-7}}.

\bibitem{simoncini2016computational}
V~Simoncini.
\newblock Computational methods for linear matrix equations.
\newblock {\em SIAM Review}, 58(3):377--441, 2016.
\newblock \href {https://doi.org/10.1137/130912839}
  {\path{doi:10.1137/130912839}}.

\bibitem{turing}
A~M Turing.
\newblock The chemical basis of morphogenesis.
\newblock {\em Bulletin of Mathematical Biology}, 52(1):153--197, 1990.
\newblock \href {https://doi.org/10.1093/oso/9780198250791.003.0022}
  {\path{doi:10.1093/oso/9780198250791.003.0022}}.

\bibitem{vanag2004waves}
V~K Vanag.
\newblock Waves and patterns in reaction--diffusion systems.
  {Belousov--Zhabotinsky} reaction in water-in-oil microemulsions.
\newblock {\em Physics-Uspekhi}, 47(9):923, 2004.
\newblock \href {https://doi.org/10.1070/pu2004v047n09abeh001742}
  {\path{doi:10.1070/pu2004v047n09abeh001742}}.

\end{thebibliography}

\end{document}